\documentclass[reqno]{amsart}

\usepackage{amssymb}
\usepackage{amsxtra}
\usepackage[all]{xy}
\usepackage{graphicx}
\usepackage{enumerate}
\usepackage{xspace}

\usepackage{hyperref}

\newtheorem{theorem}{Theorem}[section]

\newtheorem{corollary}[theorem]{Corollary}
\newtheorem{lemma}[theorem]{Lemma}
\newtheorem{proposition}[theorem]{Proposition}
\newtheorem{question}{Question}

\newtheorem{step}{Step}

\newtheorem*{lemmawon}{Lemma}
\newtheorem{maintheorem}{Main Theorem}[section]

\newtheorem*{CS}{Castelnuovo-Severi~inequality}
\newtheorem*{maroni}{Results~of~Maroni}

\theoremstyle{definition}
\newtheorem{definition}[theorem]{Definition}
\newtheorem*{definitionwon}{Definition}

\theoremstyle{remark}
\newtheorem{remark}[theorem]{Remark}
\newtheorem{remarks}[theorem]{Remarks}
\newtheorem{remarkof}{Remark}[theorem]

\numberwithin{equation}{section}

\def\sheaf#1{\ensuremath \mathcal#1}
\def\ruled#1{\ensuremath \mathbf{P}(\sheaf#1)}

\DeclareMathOperator{\bp}{Bp}
\DeclareMathOperator{\divisorgroup}{Div}%
\DeclareMathOperator{\divisor}{Divisor~of}%
\DeclareMathOperator{\gonality}{gon}
\DeclareMathOperator{\length}{length}%
\DeclareMathOperator{\pic}{Pic}%

\DeclareMathOperator{\supp}{Supp}%

\newcommand{\elm}{\ensuremath elm}%
\newcommand{\imageofblup}[1]{\ensuremath \widetilde{#1}}

\newcommand{\linsys}[1]{\ensuremath \lvert #1 \rvert}
\newcommand{\blup}{\ensuremath \varphi}
\newcommand{\bldown}{\ensuremath \psi}
\newcommand{\abs}[1]{\ensuremath \lvert #1 \rvert}
\newcommand{\PP}{\ensuremath \mathbf{P}}
\newcommand{\lineqv}{\ensuremath \sim}
\newcommand{\numeqv}{\ensuremath \equiv}

\newcommand{\nelmtrans}{\ensuremath \beta}%
\newcommand{\nelmresol}{\ensuremath \alpha}

\begin{document}
\title[Base-point-free Pencils on Triple Covers]{Base-Point-Free Pencils on Triple Covers of Smooth Curves}

\author{Dongsoo Shin}

\address{Department of Mathematical Sciences, Seoul National University, Seoul 151-747, Korea}

\email{dsshin@math.snu.ac.kr}



\date{July 25, 2007}

\begin{abstract}
Let $X$ be a smooth algebraic curve. Suppose that there exists a triple covering $f : X \to Y$ where $Y$ is a smooth algebraic curve. In this paper, we investigate the existence of morphisms from $X$ to the projective line $\mathbf{P}^1$ which do not factor through the covering $f$. For this purpose, we generalize the classical results of Maroni concerning base-point-free pencils on trigonal curves to the case of triple covers of arbitrary smooth irrational curves.
\end{abstract}

\maketitle

\section{Introduction}\label{section:Introduction}

Let $X$ be a smooth algebraic curve of genus $g_x$. Suppose that there exists a covering $f : X \to Y$ of degree $k$ where $Y$ is a smooth algebraic curve of genus $g_y$. What kind of morphisms from the cover $X$ to the projective line $\PP^1$ exist? Clearly, there are morphisms $X \to \PP^1$ induced from the base curve $Y$, that is, morphisms of the form $h' \circ f : X \to \PP^1$ for some $h' : Y \to \PP^1$. Furthermore, by the following application of the Castelnuovo-Severi inequality, every morphism $h : X \to \PP^1$ whose degree is small enough compared with the genus of $X$ is induced from the base curve $Y$ if $\deg{f}$ is a prime number.

\begin{CS}[{cf.~\cite[p.\,366]{ACGH-1985}}]
Suppose that $k$ is a prime number. Let $h : X \to \PP^1$ be a morphism of degree $d$. If
\begin{equation*}
d < \frac{g_x - kg_y + k-1}{k-1},
\end{equation*}
then $h$ factors through the covering $f$.
\end{CS}

In view of the Castelnuovo-Severi inequality, it is natural to raise the following question.

\begin{definitionwon}
A morphism $h : X \to \PP^1$ is said to be \emph{nontrivial} if
there is no morphism $h' : Y \to \PP^1$ such that $h = h' \circ f$.
\end{definitionwon}

\begin{question}\label{question:converse-of-CS-inequality}
Suppose that $k$ is a prime number. For every integer $d$ with
\begin{equation*}
d \ge \frac{g_x - kg_y + k - 1}{k-1},
\end{equation*}
does there exist a nontrivial morphism $h : X \to \PP^1$ of degree $d$?
\end{question}

For double covering case ($k=2$), Question~\ref{question:converse-of-CS-inequality} has been answered affirmatively in \cite{Keem-Ohbuchi-2004}. For trigonal curve case ($k=3$, $g_y=0$), Question~\ref{question:converse-of-CS-inequality} has been completely answered by the classical results of Maroni; cf.~\cite{Martens-Schreyer-1986}. However, if $k \ge 3$ and $g_y > 0$, then one cannot expect an affirmative answer for all that $d$'s: In the case of an affirmative answer we at least must have
\begin{equation}\label{equation:obstruction}
\frac{g_x - kg_y + k-1}{k-1} > \gonality{X} - 1,
\end{equation}
where $\gonality{X}$ is the gonality of $X$, and the bound~\eqref{equation:obstruction} is a restriction on (the gonality of) $X$ unless $k=2$ or $(k, g_y) = (3,0)$.

On the other hand, if $\frac{g_x - kg_y + k-1}{k-1} > k \cdot \frac{g_y + 3}{2} - 1$, then the inequality~\eqref{equation:obstruction} is no longer obstructional since $\gonality{X} \le k \gonality{Y} \le k \cdot \frac{g_y+3}{2}$. Thus it seems senseful to ask Question~\ref{question:converse-of-CS-inequality}, with $k \ge 3$, only for $g_x$ which is sufficiently large compared with $g_y$. For instance, it has been known that, if $X$ is a triple cover of an irrational curve $Y$ ($k=3$, $g_y > 0$) and $g_x \ge (2[(3g_y+1)/2]+1)([(3g_y+1)/2]+1)$, then there exist nontrivial morphisms of degree $d$ for all integers $d$ with $d \ge g_x - [({3g_y+1})/{2}]$; \cite[Theorem~A]{Kato-Keem-Ohbuchi-1997}. However the degree bound $d \ge g_x - [({3g_y+1})/{2}]$ is rather far from the plausible bound $d \ge ({g_x - 3g_y + 2})/{2}$ given in Question~\ref{question:converse-of-CS-inequality} and it is assumed in \cite[Theorem~A]{Kato-Keem-Ohbuchi-1997} that the base curve $Y$ is \emph{general} in the sense of Brill-Noether theorem.

In this paper, we investigate Question~\ref{question:converse-of-CS-inequality} for triple covers of arbitrary (not necessarily general) irrational smooth curves with a mild condition on the genera of the curves to be sure that the bound~\eqref{equation:obstruction} is no longer obstructional. The following Main Theorem~\ref{theorem:MainTheoremA} and \ref{theorem:MainTheoremB} give a partial answer to Question~\ref{question:converse-of-CS-inequality} in case of triple covers, which may be considered as generalizations of the results of Maroni to the case of triple covers of irrational smooth curves; see Remark~\ref{remark:final-remarks}. Recall that a certain rank two vector bundle, so-called \emph{Tschirnhausen module}, can be associated to every triple covering; cf.~\S\ref{subsection:preliminaries}.

\begin{maintheorem}[{Theorem~\ref{theorem:d>=b_2+4g_y}}]\label{theorem:MainTheoremA}%
Let $X$ be a smooth irreducible curve of genus $g_x$. Suppose that there exists a triple covering $f : X \to Y$ where $Y$ is a smooth irreducible curve of genus $g_y$ and suppose that $g_x \ge 9g_y + 4$. Let $\sheaf{E}^{\spcheck}$ be the Tschirnhausen module for the triple covering $f$ and $e$ the $e$-invariant of the ruled surface $\ruled{E}$. For every integer $d$ with
\begin{equation*}
d \ge \frac{g_x - 3g_y + 2}{2} + \frac{\abs{e}}{2} + 4g_y,
\end{equation*}
there exists a nontrivial morphism $h : X \to \PP^1$ of degree $d$.
\end{maintheorem}

\begin{maintheorem}[{Theorem~\ref{theorem:CS-inequality-for-triple-coverings}}]\label{theorem:MainTheoremB}%
Let $X$ be a smooth irreducible curve of genus $g_x$. Suppose that there exists a triple covering $f : X \to Y$ where $Y$ is a smooth irreducible curve of genus $g_y$. Let $\sheaf{E}^{\spcheck}$ be the Tschirnhausen module for the triple covering $f$ and $e$ the $e$-invariant of the ruled surface $\ruled{E}$. If there exists a nontrivial morphism $h : X \to \PP^1$ of degree $d$, then
\begin{equation*}
d \ge \frac{g_x - 3g_y + 2}{2} + \frac{\abs{e}}{2}.
\end{equation*}
\end{maintheorem}

Note that Question~\ref{question:converse-of-CS-inequality} cannot be answered affirmatively for triple coverings with $e \neq 0$ by Main Theorem~\ref{theorem:MainTheoremB}; hence we raise a weaker question as follows.

\begin{question}\label{question:weakened-version-of-converse-of-CS-inequality}
Let $Y$ be an irreducible smooth curve of genus $g_y$. Let $k \ge 2$ be a prime number and $g_x$ an integer with $g_x - kg_y + k - 1 > 0$. For every integer $d$ with
\begin{equation*}
d \ge \frac{g_x - kg_y + k - 1}{k-1},
\end{equation*}
does there exist an irreducible smooth curve $X$ of genus $g_x$ together with a degree $k$ covering $f : X \to Y$ such that $X$ admits a nontrivial morphism $h : X \to \PP^1$ of degree $d$?
\end{question}

It is clear that we must have further conditions on $g_x$ provided that Question B has an affirmative answer: Since $d \ge \gonality(X) \ge \gonality(Y)$, we have $\frac{g_x - kg_y + k-1}{k-1} > \gonality(Y)-1$, and since $\gonality(Y) \ge \frac{g_y + 2}{2}$ is possible, one obtains, for an affirmative answer for any $Y$, a bound on $g_x$ which is worse than that given by the condition $g_x - kg_y + k-1>0$.

For double covering case ($k=2$), \cite{Ballico-Keem-Park-2004} gives an affirmative answer to Question~\ref{question:weakened-version-of-converse-of-CS-inequality} for all $g_x \ge 4g_y + 5$ and $d \ge g_x - 2g_y + 1$. For trigonal curve case $(k, g_y) = (3, 0)$, it has been known that there exist trigonal curves together with nontrivial morphisms to $\PP^1$ of degree $d$ for all integer $d \ge ({g_x + 2})/{2}$ if $g_x \ge 5$; cf.~\cite{Martens-Schreyer-1986}. In case of $k = 3$ and $g_y > 0$, in \cite[Example~3.3]{Kato-Keem-Ohbuchi-1997}, the authors proved that, if $g_x \ge 7g_y - 4$, then, for every integer $d$ with
\begin{equation*}
\text{$\frac{g_x - 3g_y + 2}{2} \le d \le (g_x - 3g_y + 2)$ and $d
\equiv 2g_x - 2 \mod{3}$},
\end{equation*}
there exists an irreducible smooth curve $X$ of genus $g_x$ and a cyclic triple covering $f : X \to Y$ such that $X$ admits a nontrivial morphism $X \to \PP^1$ of degree $d$; however, the range of the degrees $d$ is somewhat restricted.

The following theorem provides a partial answer to Question~\ref{question:weakened-version-of-converse-of-CS-inequality} for triple covers of irrational smooth curves, which may be considered as a generalization of trigonal curve case.

\begin{maintheorem}[{Theorem~\ref{theorem:existence-irrational-curve}}]\label{theorem:MainTheoremC}%
Let $Y$ be an irreducible smooth curve of genus $g_y \ge 1$, and let $g_x$ be an integer with $g_x \ge 37g_y - 2$. For every integer $d$ with
\begin{equation*}
d \ge \frac{g_x - 3g_y + 2}{2} + g_y,
\end{equation*}
there exists an irreducible smooth curve $X$ of genus $g_x$ together with a triple covering $f : X \to Y$ such that $X$ admits a nontrivial morphism $h : X \to \PP^1$ of degree $d$. If $g_y \ge 5$, we get sharper lower bound for $d$:
\begin{equation*}
d \ge \frac{g_x - 3g_y + 2}{2} + \frac{g_y + 3}{2}.
\end{equation*}
\end{maintheorem}

The paper is organized as follows. In Section~\ref{section:base-point-free-pencils-on-triple-civers}, we will prove Main Theorem~\ref{theorem:MainTheoremA} and Main Theorem~\ref{theorem:MainTheoremB}. In Section~\ref{section:ruled-surface-under-elm}, we will investigate what conditions may be imposed on the Tschirnhausen module for a triple cover $f : X \to Y$ if $X$ admits a nontrivial morphism $X \to \PP^1$ of minimal possible degree given by Main Theorem~\ref{theorem:MainTheoremB}: If a triple cover $X$ admits a nontrivial morphism $X \to \PP^1$ of minimal possible degree, then its Tschirnhausen module must be decomposable; Proposition~\ref{proposition:minimal-possible=>decomposible}. Finally, in Section~\ref{section:existence-of-triple-covers}, we will prove Main Theorem~\ref{theorem:MainTheoremC}.

\subsubsection{Ideas of proofs.} Let $X$ and $Y$ be irreducible smooth curves. Assume that there exists a triple covering $f : X \to Y$. Let $\sheaf{E}^{\spcheck}$ be the Tschirnhausen module for the triple covering $f$. It has been known that there exists an embedding $i : X \hookrightarrow \ruled{E}$ of $X$ into the ruled surface $\ruled{E}$; \cite[Theorem~2.1]{Casnati-Ekedahl-1996}.

For proving Main Theorem~\ref{theorem:MainTheoremA}, we first choose a linear series $\sheaf{L}$ on $\ruled{E}$ so that $\sheaf{L}$ is very ample or separates points in the same fiber of the triple covering $f$. Then $\sheaf{L}$ cuts out a base-point-free linear series on $X$ that is not induced from the base curve $Y$. The idea is based on the proof of \cite[Lemma~1]{Martens-Schreyer-1986}.

The main tools for proving Main Theorem~\ref{theorem:MainTheoremB} are elementary transformations. It has been known that the ruled surface $\ruled{E}$ (in fact, every ruled surface) can be transformed to $Y \times \PP^1$ by a finite sequence of elementary transformations; cf.~\cite[V, Ex. 5.5]{Hartshorne-1977}. Let $\elm : \ruled{E} \dashrightarrow Y \times \PP^1$ be a sequence of elementary transformations that transforms $\ruled{E}$ to $Y \times \PP^1$. Then it is clear that the morphism
\begin{equation*}
h = p_2 \circ \elm|_X : X \to \PP^1
\end{equation*}
is nontrivial, where $p_2 : Y \times \PP^1 \to \PP^1$ is the second projection. The main ingredient of the proof of Main Theorem~\ref{theorem:MainTheoremB} is Proposition~\ref{proposition:general-form-of-bpf-pencils}: Every nontrivial morphism $h : X \to \PP^1$ is, in fact, of the form $p_2 \circ \elm|_X$ for some sequence of elementary transformations $\elm : \ruled{E} \dashrightarrow Y \times \PP^1$ that transforms $\ruled{E}$ to $Y \times \PP^1$. We prove Main Theorem~\ref{theorem:MainTheoremB} by investigating behaviors of ruled surfaces and their $e$-invariants under elementary transformations.

We prove Main Theorem~\ref{theorem:MainTheoremC} as follows. It has been known that every triple covering $f : X \to Y$ of a given curve $Y$ is determined by a rank two vector bundle $\sheaf{E}$ on $Y$ and a smooth zero locus $X$ of a section in $H^0(\ruled{E}, \pi^{\ast}{\det{\sheaf{E}}}^{-1}(3))$, where $\pi : \ruled{E} \to Y$ is the projection; \cite[Theorem~3.4]{Casnati-Ekedahl-1996}. We choose a rank two vector bundle $\sheaf{E}$ on the base curve $Y$ which is decomposable and we take some points $P_1, \dotsc, P_e \in \ruled{E}$ so that the elementary transformation $\elm$ with $P_1, \dotsc, P_e$ as centers transform $\ruled{E}$ into $Y \times \PP^1$. We then prove that there exists a irreducible smooth zero locus $X$ of a section in $H^0(\ruled{E}, \pi^{\ast}{\det{\sheaf{E}}}^{-1}(3))$ passing through the given points $P_1, \dotsc, P_e$ and we compute the degree of the morphism $h = p_2 \circ \elm|_X : X \to \PP^1$ by investigating behaviors of trisections of ruled surfaces under elementary transformations.

\subsection{Preliminaries}\label{subsection:preliminaries}

We collect some results concerning triple coverings, especially results related to Tschirnhausen modules; cf.~\cite{Casnati-Ekedahl-1996} and \cite{Miranda-1985} for detail. We briefly review basics of sections, especially minimal degree sections, of ruled surfaces; cf.~\cite{Maruyama-1970}. Finally, we recall the definition of elementary transformation; cf.~\cite[V,5.7.1]{Hartshorne-1977} and \cite{Maruyama-1970}.

\subsubsection{Triple coverings.}

Let $X$ and $Y$ be irreducible smooth curves of genus $g_x$ and $g_y$. Assume that there exists a triple covering $f : X \to Y$. One can associate a split exact sequence
\begin{equation*}
0 \to \sheaf{O_Y} \xrightarrow{f^{\sharp}} f_{\ast}\sheaf{O_X} \to
\sheaf{E}^{\spcheck} \to 0,
\end{equation*}
where $\sheaf{E}^{\spcheck}$ is a locally free $\sheaf{O_Y}$-sheaf of rank $2$ called, according to R. Miranda \cite{Miranda-1985}, \emph{the Tschirnhausen module} for the triple covering $f$; cf.~\cite{Casnati-Ekedahl-1996}, \cite{Miranda-1985}.

The branch locus of the triple covering $f$ is a divisor whose associated line bundle is $(\det{\sheaf{E}})^{\otimes 2}$; \cite[Proposition~4.7]{Miranda-1985}. Setting
\begin{equation}\label{preliminaries-equation:O_Y(B)=det(E)}
\sheaf{O_Y}(B) = \det{\sheaf{E}},
\end{equation}
we have
\begin{equation}\label{preliminaries-equation:b=g_x-3g_y+2}
b := \deg{B} = \deg{\sheaf{E}} = g_x - 3g_y + 2
\end{equation}
by the Riemann-Hurwitz formula.

There exists an embedding $i : X \hookrightarrow \ruled{E}$ of $X$ into the ruled surface $\ruled{E}$ such that $f = \pi \circ i$, where $\pi : \ruled{E} \to Y$ is the projection; cf.~\cite{Casnati-Ekedahl-1996}. Furthermore, such an embedding is unique.
\begin{lemma}[{\cite[Theorem~2.1(i)]{Casnati-Ekedahl-1996}}]\label{preliminaries-lemma:uniqueness-of-embedding}
If there is an embedding $j : X \hookrightarrow \PP$ into another ruled surface $\pi' : \PP \to Y$ such that $f = \pi' \circ j$, then $\PP \cong \ruled{E}$.
\end{lemma}

\subsubsection{Sections of ruled surfaces.}

Let $\sheaf{F}$ be any rank two vector bundle on a curve $Y$ and let $\pi : \ruled{F} \to Y$ be the projection. Suppose $\sheaf{F_0} = \sheaf{F} \otimes \sheaf{O_Y}(-N)$ is normalized. A section $S_0$ of $\ruled{F}$ is called a \emph{minimal degree section} of $\ruled{F}$ if $S_0$ is a section whose self intersection number is minimum among all sections on $\ruled{F}$. Let $e$ be the $e$-invariant of $\ruled{E}$. Since $\det{\sheaf{F_0}} = \det{\sheaf{F}} \otimes \sheaf{O_Y}(-2N)$, we have
\begin{equation}\label{preliminaries-equation:-e=b-2n}
e = - \deg{\sheaf{F_0}} = 2 \deg{N} - \deg{B} = 2n-b = - S_0^2.
\end{equation}

Giving a section of $\ruled{F}$ is equivalent to giving a subbundle of $\sheaf{F}$; cf.~\cite[Lemma~1.14]{Maruyama-1970}. More precisely, if $S$ is a section of $\ruled{F}$ and $\sheaf{L_S}$ is the corresponding subbundle, then we have an exact sequence
\begin{equation*}
0 \to \sheaf{L_S} \to \sheaf{F} \to \sheaf{L} \to 0
\end{equation*}
for some subbundle $\sheaf{L}$. The subbundle $\sheaf{L_{S_0}}$ corresponding to the minimal degree section $S_0$ is a maximal degree subbundle of $\sheaf{F}$, and vice versa; cf.~\cite[Theorem~1.16]{Maruyama-1970}. Especially, we have the following lemma.

\begin{lemmawon}[{\cite[Theorem~1.16]{Maruyama-1970}}]
Let $S$ be a section of $\ruled{F}$ and $\sheaf{L_S}$ the subbundle of $\sheaf{E}$ associated to the section $S$. Then we have
\begin{equation}\label{preliminaries-equation:pi(S.S)}
\pi_{\ast}(S \cdot S) = \divisor{(\det{\sheaf{E}} \otimes \sheaf{L_S}^{-2})},
\end{equation}
where $(S \cdot S)$ denotes the intersection divisor on the smooth curve $S$.
\end{lemmawon}

If $\sheaf{F}$ is decomposable, then there exist two disjoint sections; cf.~\cite[V,Ex.2.2]{Hartshorne-1977}. We record this fact for the convenience of the reader.

\begin{lemma}\label{lemma:two-disjoint-sections}
Let $\sheaf{F} = \sheaf{O_Y}(B_1) \oplus \sheaf{O_Y}(B_2)$ be a decomposable rank two vector bundle on a smooth curve $Y$ with $\deg{B_1} \le \deg{B_2}$. Let $S_0$ be a minimal degree section of $\ruled{F}$ and $\pi : \ruled{F} \to Y$ the projection. Then there exists a section $S \in |S_0 + \pi^{\ast}(B_2 - B_1)|$ such that $S_0 \cap S = \varnothing$. Conversely, if $S$ is a section such that $S \cap S_0 = \varnothing$, then $S \in |S_0 + \pi^{\ast}(B_2 - B_1)|$.
\end{lemma}

\begin{proof}
Let $S$ be a section corresponding to the subbundle $\sheaf{O_Y}(B_1)$ of $\sheaf{F}$. Consider $(S \cdot S_0)$ as a divisor on the section $S_0$. According to \cite[Remark~1.19]{Maruyama-1970} and by \eqref{preliminaries-equation:pi(S.S)}, we have
\begin{equation*}
2\pi_{\ast}(S \cdot S_0) \lineqv \pi_{\ast}(S \cdot S) + \pi_{\ast}(S_0 \cdot S_0)  \lineqv (B_2 - B_1) + (B_1 - B_2) \lineqv 0;
\end{equation*}
hence $S \cap S_0 = \varnothing$. Conversely, suppose that $S$ is a section such that $S_0 \cap S = \varnothing$. Then $\pi_{\ast}(S \cdot S_0) \lineqv 0$. By \eqref{preliminaries-equation:pi(S.S)}, we have
\begin{equation*}
\pi_{\ast}(S_0 \cdot S_0) = \divisor{(\det{\sheaf{F}} \otimes \sheaf{L_S}^{-2})} \lineqv B_1 - B_2.
\end{equation*}
Set $S \lineqv S_0 + \pi^{\ast}D$ for some divisor $D$ on $Y$. Then we have
\begin{equation*}
0 \lineqv \pi_{\ast}(S \cdot S_0) \lineqv \pi_{\ast}((S_0 + \pi^{\ast}{D}) \cdot S_0 )\lineqv (B_1 - B_2) + D;
\end{equation*}
therefore $D \lineqv B_2 - B_1$; hence $S \in |S_0 + \pi^{\ast}(B_2 - B_1)|$.
\end{proof}

Return to the triple covering case. Suppose that
\begin{equation}\label{preliminaries-equation:E_0=E-times-O_Y(-N)}
\sheaf{E_0} = \sheaf{E} \otimes \sheaf{O_Y}(-N)
\end{equation}
is normalized; cf.~\cite[V,~2.8]{Hartshorne-1977}. For a canonical divisor $K_{\ruled{E}}$ of $\ruled{E}$ and for the triple cover $X$ regarded as a divisor of $\ruled{E}$, we have
\begin{equation}\label{preliminaries-equation:K_E-and-X}
\text{$K_{\ruled{E}} \lineqv -2S_0 + \pi^{\ast}(B-2N + K_Y)$ and $X \lineqv 3S_0 + \pi^{\ast}(3N-B)$},
\end{equation}
where $K_Y$ is a canonical divisor of $Y$; cf.~\cite[V,~2.10]{Hartshorne-1977}, \cite[Theorem~2.1]{Casnati-Ekedahl-1996}.

\subsubsection{Elementary transformations.}

Let $W$ be a ruled surface with the projection $\pi : W \to Y$. Let $P \in W$ and $F=\pi^{-1}(\pi(P))$. Let $\blup : \imageofblup{W} \to W$ be the blowing up of $W$ at $P$. The strict transform $\imageofblup{F}$ of the fiber $F$ by the blowing up $\blup$ on $\imageofblup{W}$ is isomorphic to $\PP^1$ and $\imageofblup{F}^2 = -1$. Therefore we can blow down $\imageofblup{W}$ along $\imageofblup{F}$; in other words, there is a morphism $\bldown : \imageofblup{W} \to W'$ to a surface $W'$ which is the blowing up of $W'$ contracting $\imageofblup{F}$ to a point $P'$. It is easy to prove that the surface $W'$ is again a ruled surface over $Y$. The new ruled surface $W'$ is called \emph{the elementary transformation of $W$ with center $P$}. We denote the birational map $\bldown \circ \blup^{-1}$ by $\elm_P$ and call it \emph{the elementary transform with center $P$}. The surface $W'$ is denoted by $\elm_P(W)$ and, if no confusion is likely to occur, we denote the projection $\elm_{P}(W) \to Y$ by $\pi$ again. Note that $\elm_{P'}$ is the inverse of $\elm_P$.

Let $P_1, \dotsc, P_n$ be distinct points in $W$ such that $\pi(P_i) \neq \pi(P_j)$ for $i \neq j$. We inductively define $\elm_{P_1, \dotsc, P_n}$ by
\begin{equation*}
\elm_{P_1, \dotsc, P_n} = \elm_{P_n} \circ \elm_{P_1, \dotsc,
P_{n-1}},
\end{equation*}
where $P_{n}$ is considered as the point $\elm_{P_1, \dotsc, P_{n-1}}(P_n) \in \elm_{P_1, \dotsc, P_{n-1}}(W)$. Note that $\elm_{P, Q} = \elm_{Q, P}$ if $\pi(P) \neq \pi(Q)$. Let $S$ be a section of the ruled surface $W$. Then the image $\elm_{P_1, \dotsc, P_n}(S)$ of $S$ is a section of $\elm_{P_1, \dotsc, P_n}(W)$. Suppose $P \in S$. We define $\elm_{nP}$ for $n \in \mathbb{N}$ by
\begin{equation*}
\elm_{nP} = \elm_{P_1, P_2, \dotsc, P_n},
\end{equation*}
where $P_1 = P$ and $P_{m+1}$ ($m=1, \dotsc, n-1$) is the unique point in $\elm_{mP}(S) \cap \pi^{-1}(\pi(P))$. It is clear how to define $\elm_{D}$ for an  arbitrary effective divisor $D$ on a section $S$ of $W$.

\subsubsection{Notations and conventions.}

We adopt and use almost all the notations and conventions in \cite{ACGH-1985} and \cite{Hartshorne-1977}. If no confusion is likely to occur, we denote $H^{i}(C, \sheaf{O_C}(D))$ by $H^{i}(C, D)$ and $\dim{H^{i}(C, \sheaf{O_C}(D))}$ by $h^i(C, D)$. For a linear series $\linsys{D}$, $\bp(\linsys{D})$ denotes the set of all base points of $\linsys{D}$; possibly $\bp(\linsys{D}) = \varnothing$.

If we denote $\elm_{P} = \bldown \circ \blup^{-1}$ for $\elm_{P}$, then the map $\blup : \imageofblup{W} \to W$ always denotes the blowing up of $W$ at $P$ and the map $\bldown : \imageofblup{W} \to W'$ always denotes the blowing down of $\imageofblup{W}$ along the strict transform $\imageofblup{F}$ under $\blup$ of the fiber $F = \pi^{-1}(P)$. For a curve $C \subset W$, the curve $\imageofblup{C} \subset \imageofblup{W}$ denotes the \emph{strict transform} of $C$ by the blowing up $\blup$.

The integer $r_P$ is the \emph{multiplicity} of $C$ at $P$ and the integer $\delta_{P}$ is the \emph{measure of a curve singularity}, that is, $\delta_{P} = \length(\widetilde{\sheaf{O}}_{C,P}/\sheaf{O_{C,P}})$, where$\widetilde{\sheaf{O}}_{C,P}$ is the integral closure of $\sheaf{O_{C,P}}$.

\section{Base-point-free pencils on triple covers}\label{section:base-point-free-pencils-on-triple-civers}

Let $X$ and $Y$ be irreducible smooth curves of genus $g_x$ and $g_y$ and assume that there exists a triple covering $f : X \to Y$. Let $\sheaf{E}^{\spcheck}$ be the Tschirnhausen module for the triple covering $f$, $e$ the $e$-invariant of $\ruled{E}$, and $\pi : \ruled{E} \to Y$ the projection. Finally, let $p_1 : Y \times \PP^1 \to Y$ and $p_2 : Y \times \PP^1 \to \PP^1$ be the projections.

For convenience, we define the following invariant which will be used several times; cf.~Remark~\ref{remark:final-remarks}.

\begin{definition}\label{definition:M-invariant}
The \emph{M-invariant} $m$ of the triple covering $f$ is defined by the number
\begin{equation*}
m = \frac{g_x - 3g_y + 2}{2} - \frac{\abs{e}}{2} - 2.
\end{equation*}
\end{definition}

\begin{remarkof}\label{remark:Maroni-invariant}
The number $m$ is an integer since, by \eqref{preliminaries-equation:b=g_x-3g_y+2} and \eqref{preliminaries-equation:-e=b-2n}, we have
\begin{equation}\label{equation:m-is-an-integer}
m = \deg{N} - \frac{e + \abs{e}}{2} - 2.
\end{equation}
\end{remarkof}

\subsection{Existence of nontrivial morphisms}\label{subsection:existence}

We will prove the existence of nontrivial morphisms, which give a partial answer to Question~\ref{question:converse-of-CS-inequality}; Theorem~\ref{theorem:d>=b_2+4g_y}. We need the following lemma.

\begin{lemma}
Let $e$ be the $e$-invariant of $\ruled{E}$. Then
\begin{equation}\label{preliminaries-equation:|e|<=max(2g_y,(g_x-3g_y+2)/3)}
- g_y \le e \le \max \left\{2g_y - 2, \frac{g_x - 3g_y + 2}{3} \right\}.
\end{equation}
\end{lemma}

\begin{proof}
Suppose that the Tschirnhausen module $\sheaf{E}^{\spcheck}$ is indecomposable; then the $e$-invariant $e$ (in fact, any $e$-invariants of indecomposable ruled surfaces) satisfies the inequalities
\begin{equation}\label{preliminaries-equation:-g_y<=e<=2g_y-2}
- g_y \le e \le 2g_y - 2
\end{equation}
by Nagata's Theorem; cf.~\cite[V,~Ex.2.5]{Hartshorne-1977}. Suppose that the Tschirnhausen module $\sheaf{E}^{\spcheck}$ is decomposable with $\sheaf{E} = \sheaf{L} \oplus \sheaf{M}$ and $\deg{\sheaf{L}} \le \deg{\sheaf{M}}$. Since $(\det{\sheaf{E}})^{\otimes 2} = \sheaf{L}^{2} \otimes \sheaf{M}^{2}$, it follows that
\begin{equation*}
\deg{\sheaf{L}} + \deg{\sheaf{M}} = g_x - 3g_y + 2
\end{equation*}
by \eqref{preliminaries-equation:b=g_x-3g_y+2}. We also have $2\deg{\sheaf{L}} \ge \deg{\sheaf{M}}$ and $2\deg{\sheaf{M}} \ge \deg{\sheaf{L}}$ by \cite[P.1145]{Miranda-1985}. Therefore we have
\begin{equation}\label{preliminaries-equation:2L>=M-and-2M>=L}
\frac{g_x - 3g_y + 2}{3} \le \deg{\sheaf{L}} \le \frac{g_x - 3g_y + 2}{2} \le \deg{\sheaf{M}} \le \frac{2(g_x - 3g_y + 2)}{3}.
\end{equation}
Since $e = \deg{\sheaf{M}} - \deg{\sheaf{L}}$, it follows that
\begin{equation}\label{preliminaries-equation:e<=(g_x-3g_y+2)/3}
0 \le e \le \frac{g_x-3g_y+2}{3}.
\end{equation}
Hence \eqref{preliminaries-equation:|e|<=max(2g_y,(g_x-3g_y+2)/3)} follows from \eqref{preliminaries-equation:-g_y<=e<=2g_y-2} and \eqref{preliminaries-equation:e<=(g_x-3g_y+2)/3}.
\end{proof}

To get base-point-free pencils on $X$, we first take certain linear series on $\ruled{E}$.

\begin{lemma}\label{lemma:cutout-by-large-linear-series}
Assume that $g_x \ge 9g_y + 4$. Let $A \in \divisorgroup(Y)$ be an effective divisor with $\deg{A} \ge 2g_y - 1$ and let $m$ be the M-invariant of the triple covering $f$. If $\deg{A} \le m$, then the linear series $\linsys{K_X - f^{\ast}{A}}$ on $X$ is cut out by the linear series $\linsys{K_{\ruled{E}} + X - \pi^{\ast}{A}}$ on $\ruled{E}$.
\end{lemma}

\begin{remarkof}
By \eqref{preliminaries-equation:|e|<=max(2g_y,(g_x-3g_y+2)/3)}, we have $e \le ({g_x - 3g_y + 2})/{3}$; hence it follows by the assumption $g_x \ge 9g_y + 4$ that $m \ge 2g_y$. Therefore there exists a divisor $A \in \divisorgroup(Y)$ with $2g_y - 1 \le \deg{A} \le m$.
\end{remarkof}

\begin{proof}[Proof of Lemma~\ref{lemma:cutout-by-large-linear-series}.]
Let $S_0$ be a minimal degree section of $\ruled{E}$ and set $a = \deg{A}$. To begin with, we count the dimension of $\linsys{K_{\ruled{E}} + X - \pi^{\ast}{A}}$. By \eqref{preliminaries-equation:K_E-and-X}, we have
\begin{equation}\label{equation:K_E-and-K_X}
K_{\ruled{E}} + X - \pi^{\ast}{A} \lineqv S_0 + \pi^{\ast}(N + K_Y -
A).
\end{equation}
By \cite[V,~2.4]{Hartshorne-1977} and \cite[II,~7.11]{Hartshorne-1977}, we have
\begin{equation*}
h^0(\ruled{E}, \sheaf{O_{\ruled{E}}}(S_0)) = h^0(Y, \pi_{\ast}{\sheaf{O_{\ruled{E}}}(S_0)}) = h^0(Y, \sheaf{E_0}).
\end{equation*}
It follows by the projection formula that
\begin{equation*}
\begin{split}
h^0(\ruled{E}, K_{\ruled{E}} + X - \pi^{\ast}{A}) &= h^0(Y, \pi_{\ast}(\sheaf{O_{\ruled{E}}}(S_0) \otimes \pi^{\ast}\sheaf{O_{\ruled{E}}}(N + K_Y - A))) \\ %
&= h^0(Y, \sheaf{E_0} \otimes \sheaf{O_Y}(N + K_Y - A)).
\end{split}
\end{equation*}
Since $\deg(\sheaf{E_0} \otimes \sheaf{O_Y}(N + K_Y - A)) = -e + 2(n
+ 2g_y - 2 - a)$, it follows by Riemann-Roch Theorem,
\eqref{preliminaries-equation:b=g_x-3g_y+2} and
\eqref{preliminaries-equation:-e=b-2n}, and Serre duality that
\begin{multline}\label{equation:h0(K_E+X-piA)-first}
h^0(\ruled{E}, K_{\ruled{E}} + X - \pi^{\ast}{A}) \\ %
\begin{aligned}
&= -e + 2(n + 2g_y - 2 - a) + 2(1-g_y) + h^1(Y, \sheaf{E_0} \otimes \sheaf{O_Y}(N + K_Y - A)) \\ %
&= g_x - g_y - 2a + h^0(Y, \sheaf{E_0}^{\spcheck} \otimes
\sheaf{O_Y}(-N + A)).
\end{aligned}
\end{multline}
To compute $h^0(Y, \sheaf{E_0}^{\spcheck} \otimes \sheaf{O_Y}(-N +
A))$, we have
\begin{equation}\label{equation:E_0^spcheck-cong-E_0-otimes-blabla}
\begin{split}
\sheaf{E_0}^{\spcheck} \cong \sheaf{E_0} \otimes
\det{\sheaf{E_0}}^{-1} \cong \sheaf{E_0} \otimes
(\det{\sheaf{E}}^{-1} \otimes \sheaf{O_Y}(2N)) = \sheaf{E_0} \otimes
\sheaf{O_Y}(-B + 2N).
\end{split}
\end{equation}
Therefore it follows by
\eqref{equation:h0(K_E+X-piA)-first} and
\eqref{equation:E_0^spcheck-cong-E_0-otimes-blabla} that
\begin{multline}\label{equation:h0(K_E+X-piA)-second}
h^0(\ruled{E}, K_{\ruled{E}} + X - \pi^{\ast}{A}) \\ %
= g_x - g_y - 2a + h^0(Y, \sheaf{E_0} \otimes \sheaf{O_Y}(-B + N +
A)).
\end{multline}
By \eqref{preliminaries-equation:-e=b-2n} and
\eqref{equation:m-is-an-integer}, we have
\begin{equation*}
\deg(-B + N + A) \le -b + n + m \le -2 - \frac{\abs{e} - e}{2} <
0.
\end{equation*}
Since $\sheaf{E_0}$ is normalized, we have
\begin{equation}\label{equation:h0(E_0+-B+N+A)}
h^0(Y, \sheaf{E_0} \otimes \sheaf{O_Y}(-B + N + A)) = 0.
\end{equation}
Therefore it follows by
\eqref{equation:h0(K_E+X-piA)-second} and
\eqref{equation:h0(E_0+-B+N+A)} that
\begin{equation}\label{equation:h0(K_E+X-piA)}
h^0(\ruled{E}, K_{\ruled{E}} + X - \pi^{\ast}{A}) = g_x - g_y - 2a.
\end{equation}

We now prove that
\begin{equation*}
h^0(X, K_{X} - f^{\ast}{A}) = h^0(\ruled{E}, K_{\ruled{E}} + X -
\pi^{\ast}{A}) = g_x - g_y - 2a.
\end{equation*}
Since $f$ is finite, we have $h^0(X, \sheaf{O_X}(f^{\ast}{A})) =
h^0(Y, f_{\ast}(\sheaf{O_X}(f^{\ast}{A})))$ by \cite[III,~Ex.4.1]{Hartshorne-1977}. Hence, by the projection formula, we have
\begin{equation}\label{equation:h0(fA)}
\begin{split}
&h^0(X, \sheaf{O_X}(f^{\ast}{A})) = h^0(Y, f_{\ast}(\sheaf{O_X}(f^{\ast}{A}) \otimes \sheaf{O_X})) = h^0(Y, \sheaf{O_Y}(A) \otimes f_{\ast}{\sheaf{O_X}}) \\ %
&\qquad= h^0(Y, \sheaf{O_Y}(A) \otimes (\sheaf{O_Y} \oplus \sheaf{E}^{\spcheck})) = h^0(Y, \sheaf{O_Y}(A)) + h^0(Y, \sheaf{E}^{\spcheck} \otimes \sheaf{O_Y}(A)) \\ %
&\qquad= a - g_y + 1 + h^0(Y, \sheaf{E}^{\spcheck} \otimes
\sheaf{O_Y}(A)),
\end{split}
\end{equation}
where $h^1(Y, \sheaf{O_Y}(A)) = 0$ by the assumption $\deg{A} \ge
2g_y - 1$. It remains to count $h^0(Y, \sheaf{E}^{\spcheck} \otimes
\sheaf{O_Y}(A))$. We have $\sheaf{E_0}^{\spcheck} \cong \sheaf{E_0}
\otimes \sheaf{O_Y}(-B + 2N)$ by
\eqref{equation:E_0^spcheck-cong-E_0-otimes-blabla}, but
$\sheaf{E_0} = \sheaf{E} \otimes \sheaf{O_Y}(-N)$; hence
$\sheaf{E}^{\spcheck} \cong \sheaf{E_0} \otimes \sheaf{O_Y}(-B +
N)$. Therefore it follows that
\begin{equation*}
h^0(Y, \sheaf{E}^{\spcheck} \otimes \sheaf{O_Y}(A)) = h^0(Y,
\sheaf{E_0} \otimes \sheaf{O_Y}(-B + N + A)) = 0
\end{equation*}
by \eqref{equation:h0(E_0+-B+N+A)}.
Then we have
\begin{equation*}
h^0(X, \sheaf{O_X}(f^{\ast}{A})) = a - g_y + 1
\end{equation*}
by \eqref{equation:h0(fA)}; hence
\begin{equation}\label{equation:h0(K_X-fA)}
h^0(X, K_X - f^{\ast}{A}) = g_x - g_y - 2a.
\end{equation}
Thus it follows by \eqref{equation:h0(K_E+X-piA)} and
\eqref{equation:h0(K_X-fA)} that
\begin{equation}\label{equation:h0(K_E+X-piA)=h0(K_X-fA)}
h^0(\ruled{E}, K_{\ruled{E}} + X - \pi^{\ast}{A}) = h^0(X, K_X -
f^{\ast}{A}) = g_x - g_y - 2a.
\end{equation}

Finally, we claim that the following restriction map $\gamma$ is
injective:
\begin{equation*}
\gamma : H^0(\ruled{E}, K_{\ruled{E}} + X - \pi^{\ast}{A}) \to
H^0(X, K_X - f^{\ast}{A}).
\end{equation*}
Let $C, C' \in \linsys{K_{\ruled{E}} + X - \pi^{\ast}{A}}$. By
\eqref{equation:m-is-an-integer}, $m \le n - e - 2$; hence
\begin{equation}\label{equation:a<=m<=n-e-2<=n-e+2g_y-2}
a \le m \le n - e - 2 \le n - e + 2g_y - 2.
\end{equation}
We have $S_0^2 = b - 2n$ by
\eqref{preliminaries-equation:-e=b-2n}; it follows by
\eqref{preliminaries-equation:b=g_x-3g_y+2} and
\eqref{equation:K_E-and-K_X} that
\begin{equation*}
C.X - C.C' = 2n - e + 2g_y - 2 - a.
\end{equation*}
By \eqref{preliminaries-equation:b=g_x-3g_y+2} and
\eqref{preliminaries-equation:-e=b-2n}, we have $e = 2n - (g_x -
3g_y + 2)$, and, by
\eqref{preliminaries-equation:|e|<=max(2g_y,(g_x-3g_y+2)/3)},
we have $e \ge -g_y$. Therefore it follows that $2n - (g_x - 3g_y +
2) \ge -g_y$; hence it follows by the assumption $g_x \ge 9g_y + 4$
that $n = \deg{N} > 0$. Furthermore we have $n - e + 2g_y - 2-a \ge
0$ by \eqref{equation:a<=m<=n-e-2<=n-e+2g_y-2}; hence it
follows that
\begin{equation*}
C.X - C.C' = n + (n - e + 2g_y - 2 - a) > 0.
\end{equation*}
Therefore the restriction map $\gamma$ is injective as asserted;
hence the linear series $\linsys{K_X - f^{\ast}{A}}$ is cut out by
the linear series $\linsys{K_{\ruled{E}} + X - \pi^{\ast}{A}}$.
\end{proof}

We will get base-point-free linear series on $X$ which are not
induced from the base curve $Y$.

\begin{lemma}\label{lemma:very-ample-or-simple}
Assume that $g_x \ge 9g_y + 4$. Let $A \in \divisorgroup(Y)$ be an
effective divisor with $\deg{A} \ge 2g_y - 1$ and let $m$ be the
M-invariant of the triple covering $f$. If $\deg{A} \le m - 1$, then
the linear series $\linsys{K_X - f^{\ast}{A}}$ on $X$ is very ample.
If $\deg{A} = m$, then $\linsys{K_X - f^{\ast}{A}}$ is
base-point-free and separates points in the same fiber of the triple
covering $f$.
\end{lemma}

\begin{proof}
Set $\sheaf{L} = \linsys{K_{\ruled{E}} + X - \pi^{\ast}{A}}$. By
\eqref{preliminaries-equation:K_E-and-X}, we have
$\sheaf{L} = \linsys{S_0 + \pi^{\ast}(N + K_Y - A)}$. We first prove
that $\sheaf{L}$ is very ample if $\deg{A} \le m - 1$, or
$\sheaf{L}$ is base-point-free and separates points in the same
fiber of the triple covering $f$ if $\deg{A} = m$.

Suppose that $a = \deg{A} \le m - 1$. We have $m = n - 2
-(e+\abs{e})/2$ and $e=2n-b$ by
\eqref{preliminaries-equation:-e=b-2n} and
\eqref{equation:m-is-an-integer}; it follows that
\begin{align*}
&\deg(N + K_Y - A) \ge 2g_y + 1 + \frac{\abs{e} + e}{2} \ge 2g_y + 1,\\ %
&\deg(N + K_Y - A + (B-2N)) \ge 2g_y + 1 + \frac{\abs{e} - e}{2}
\ge 2g_y + 1.
\end{align*}
Therefore $\linsys{N + K_Y - A}$ and $\linsys{N + K_Y - A + (B-2N)}$
are very ample, and $\linsys{N + K_Y - A - P}$ and $\linsys{N + K_Y
- A + (B-2N) - P}$ are nonspecial for every $P \in Y$. Thus
$\sheaf{L}$ is very ample by \cite[V,~Ex.2.11]{Hartshorne-1977}.

Suppose that $a =\deg{A} = m$. Then we have
\begin{align*}
&\deg(N + K_Y - A) \ge 2g_y + \frac{\abs{e} + e}{2} \ge 2g_y,\\ %
&\deg(N + K_Y - A + (B-2N)) \ge 2g_y + \frac{\abs{e} - e}{2} \ge
2g_y,
\end{align*}
by \eqref{preliminaries-equation:b=g_x-3g_y+2},
\eqref{preliminaries-equation:-e=b-2n}, and
\eqref{equation:m-is-an-integer}. Therefore $\linsys{N + K_Y - A}$
and $\linsys{N + K_Y - A + (B-2N)}$ have no base points, and
$\linsys{N + K_Y - A}$ is nonspecial. Thus $\sheaf{L}$ is
base-point-free by \cite[V,~Ex.2.11]{Hartshorne-1977}. Furthermore
$\sheaf{L}$ separates points in the same fiber of the triple
covering $f$ because $\linsys{S_0 + \pi^{\ast}(N + K_Y - A)}$
contains a section by \cite[V,~Ex.2.11]{Hartshorne-1977}.

Since $\sheaf{L}$ is base-point-free for any effective divisor $A
\in \divisorgroup(Y)$ with $2g_y - 1 \le \deg{A} \le m$, the linear
series $\linsys{K_X - f^{\ast}{A}}$ on $X$ is also base-point-free
because $\linsys{K_X - f^{\ast}{A}}$ is cut out by $\sheaf{L}$ by
Lemma~\ref{lemma:cutout-by-large-linear-series}. Set
\begin{equation*}
r = \dim{\linsys{K_X - f^{\ast}{A}}} = \dim{\linsys{K_{\ruled{E}} +
X - \pi^{\ast}{A}}};
\end{equation*}
cf.~\eqref{equation:h0(K_E+X-piA)=h0(K_X-fA)}. Let $\alpha
: X \to \PP^r$ and $\beta : \ruled{E} \to \PP^r$ be the morphisms
associated to the base-point-free linear series $\linsys{K_X -
f^{\ast}{A}}$ and $\sheaf{L}$, respectively. The morphism $\beta$ is
an embedding for $a \le m - 1$ or a morphism which separates points
in the same fiber of the triple covering $f$ for $a = m$; but, we
have $\alpha = \beta \circ i$, where $i : X \hookrightarrow
\ruled{E}$ is the inclusion. Therefore the morphism $\alpha$ is also
an embedding for $a \le m - 1$ or a morphism which separates points
in the same fiber of the triple covering $f$ for $a = m$.
\end{proof}

We need the following lemma.

\begin{lemma}[{\cite[Lemma~2.2.1]{Coppens-Keem-Martens-1992}}]\label{lemma:subtracting-general-points}
Let $f : X \to Y$ be a covering of degree $k \ge 2$ of smooth curves
and let $g_{d}^{r}$($r \ge 1$) be a base-point-free linear series on
$X$ which is not induced by $Y$. Let $P_1, \dotsc, P_{r-1}$ be $r-1$
general points on $X$. Then the base-point-free part of the pencil
$g_{d}^{r}(-P_1 - \dotsb - P_{r-1})$ is not induced by $Y$.
\end{lemma}

We now prove the existence theorem.

\begin{theorem}\label{theorem:d>=b_2+4g_y}
Let $X$ be a smooth irreducible curve of genus $g_x$. Suppose that
there exists a triple covering $f : X \to Y$ where $Y$ is a smooth
irreducible curve of genus $g_y$ and suppose that $g_x \ge 9g_y +
4$. Let $\sheaf{E}^{\spcheck}$ be the Tschirnhausen module for the
triple covering $f$ and $e$ the $e$-invariant of the ruled surface
$\ruled{E}$. For every integer $d$ with
\begin{equation*}
d \ge \frac{g_x - 3g_y + 2}{2} + \frac{\abs{e}}{2} + 4g_y,
\end{equation*}
there exists a nontrivial morphism $h : X \to \PP^1$ of degree $d$.
\end{theorem}

\begin{proof}
It is obvious for $d \ge g_x + 1$; just take a general $\linsys{D}
\in \pic(X)$ which is nonspecial. From now on, assume that $d \le
g_x$. Let $m$ be the M-invariant of the triple covering $f$. Choose
an effective divisor $A \in \divisorgroup(Y)$ of degree $a$ with
$2g_y - 1 \le a \le m$. By Lemma~\ref{lemma:very-ample-or-simple},
the linear series $\linsys{K_X - f^{\ast}{A}}$ is not composed with
the triple covering $f$. By \eqref{equation:h0(K_X-fA)}, we
have
\begin{equation*}
\dim \linsys{K_X - f^{\ast}{A}} = g_x - g_y - 2a - 1.
\end{equation*}
Therefore, subtracting general $g_x - g_y - 2a - 2$ points from
$\linsys{K_X - f^{\ast}{A}}$, we have a base-point-free pencil of
degree $g_x + g_y - a$, which is not composed with the triple
covering $f$ by Lemma~\ref{lemma:subtracting-general-points}. Since
\begin{equation*}
2g_y - 1 \le a \le m = \frac{g_x - 3g_y + 2}{2} - \frac{\abs{e}}{2}
- 2,
\end{equation*}
we have
\begin{equation*}
\frac{g_x - 3g_y + 2}{2} + \frac{\abs{e}}{2} + 4g_y \le g_x + g_y -
a \le g_x - g_y + 1.
\end{equation*}
Therefore, for every integer $d$ with
\begin{equation*}
\frac{g_x - 3g_y + 2}{2} + \frac{\abs{e}}{2} + 4g_y \le d \le g_x -
g_y + 1,
\end{equation*}
there exists a base-point-free pencil of degree $d$ which is not
composed with the triple covering $f$.

If $g_y = 0$, the proof is done. Assume that $g_y \ge 1$. We need
the following lemma; here we adopt the conventions and notation used
in \cite{ACGH-1985}.

\begin{lemmawon}[{\cite[Lemma~3.2]{Ballico-Keem-Park-2004}}]
Fix an integer $s \ge 1$. Let $C$ be a smooth curve of genus $g \ge
4s-4$ defined over an algebraically closed field of characteristic
zero. For an integer $d$, let $\Sigma_{d}^{1}$ be the union of
components of $W_{d}^{1}(C)$ whose general element is
base-point-free and complete. If $\Sigma_{g - s + 1}^{1} \neq
\varnothing$, then we have $\Sigma_{g - s + 2}^{1} \neq
\varnothing$.
\end{lemmawon}

Taking $C = X$ and $s = g_y$ in the above Lemma, it follows that
$\Sigma_{g_x - g_y + 1}^{1} \neq \varnothing$; hence $\Sigma_{g_x -
g_y + 2}^{1} \neq \varnothing$. By taking $s' = g_y - 1$, we again
have $\Sigma_{g_x - s' + 2} = \Sigma_{g_x - g_y + 3} \neq
\varnothing$; note that $g_x \ge 4g_y - 4 > 4s' - 4$ by the
hypothesis $g_x \ge 9g_y + 4$. We may continue this process by
taking smaller $s$'s until $s=2$ and we stop.
\end{proof}

\begin{remarkof}\label{remark:referee}
Let $f : X \to Y$ be a covering of degree $k$. For prime $k$, assume that $g_x \ge k^2 g_y + (k-1)^2$ which for $k=3$ is the bound of the above theorem. Then it follows that
\begin{equation*}
\frac{g_x - kg_y + k-1}{k-1} \ge k(g_y + 1) \ge k \gonality(Y) \ge \gonality(X).
\end{equation*}
If $\frac{g_x - kg_y + k-1}{k-1} > \gonality(X)$, then any morphism of degree $\gonality(X)$ onto $\PP^1$ factors through $f$; in particular, $\gonality(X) = k \cdot \gonality(Y)$. On the other hand, if $\frac{g_x - kg_y + k-1}{k-1} = \gonality(X)$, then $\gonality(X) = k(g_y + 1)$; hence
\begin{equation*}
k \gonality(Y) = \gonality(X) = k(g_y + 1)
\end{equation*}
i.e., $\gonality(Y) = g_y + 1$ which is only possible for $g_y \le 1$. But for $g_y \le 1$ the above equality shows that $\gonality(X) = k \gonality(Y)$, again. Therefore, assuming $g_x \ge 9g_y + 4$ in Main Theorem~\ref{theorem:MainTheoremA}, the gonality of $X$ is fixed by that of $Y$. Furthermore, the bound implies that $Y$ is the only curve of genus at most $g_y$ triply covered by $X$.
\end{remarkof}

\subsection{Castelnuovo-Severi inequality for triple coverings}

We will improve the Castelnuovo-Severi inequality in case of triple
coverings; Theorem~\ref{theorem:CS-inequality-for-triple-coverings}.
As a by-product, we may conclude that
Question~\ref{question:converse-of-CS-inequality} does not have an
affirmative answer for certain triple coverings.

First, we prove that every nontrivial morphism $X \to \PP^1$ is
determined by a certain finite sequence of elementary
transformations;
Proposition~\ref{proposition:general-form-of-bpf-pencils}. Recall
that \emph{trisections} are irreducible curves on ruled surfaces
whose intersection number with a fiber is $3$.

\begin{remark}\label{remark:X'-is-a-trisection}
Let $h : X \to \PP^1$ be a nontrivial morphism and set $X' = (f
\times h)(X) \subset Y \times \PP^1$. Since $\deg{f} = 3$ and
$\deg{h} = d$, we have $X' \numeqv 3S + dF$, where $S = p_2^{-1}(s)$
for some $s \in \PP^1$ and $F = p_1^{-1}(y)$ for some $y \in Y$.
Therefore $X'$ is a trisection of $\ruled{E}$. Let $p_a(X')$ be the
arithmetic genus of $X'$. By the adjunction formula, we have
$p_a(X') = 2d + 3g_y - 2$. Since the morphism $h$ does not factor
through the triple covering $f$, the image $X'$ is birational to
$X$; hence $X$ is the normalization of $X'$. Therefore it follows
that
\begin{equation}\label{equation:p_a(X')=2d+3g_y-2}
p_a(X') = 2d + 3g_y - 2 = g_x + \sum_{P \in X'} \delta_{P}.
\end{equation}
\end{remark}

\begin{lemma}\label{lemma:elm_at_singular_point}
Let $C$ be a trisection on a ruled surface $W$ over $Y$ and let $\pi
: W \to Y$ be the projection. For $P \in C \subset W$, set $\elm_{P}
= \bldown \circ \blup^{-1}$. Suppose that $P$ is a singular point of
$C$ with multiplicity $r_P$ and suppose that there exists a
infinitely near singular point $Q$ of $\imageofblup{C}$ lying over
$P$. Then $Q$ is the unique singular point of $\imageofblup{C}$
among points lying over $P$, and $P' = \bldown(Q)$ is the unique
singular point of $C'$ among points on $\pi^{-1}(\pi(P')) \cap C'$.
Furthermore, if $r_P = 2$ then $r_{P'}=2$ and $\delta_{P'} =
\delta_{P} - 1$, and if $r_P = 3$ then $\delta_{P'} = \delta_{P} -
3$.
\end{lemma}

\begin{proof}
Let $E$ be the exceptional divisor of the blowing up $\blup$. Set $F
= \pi^{-1}(\pi(P))$, $W' = \elm_{P}(W)$, $C' = \elm_{P}(C)$, and $E'
= \psi(E)$. Suppose that $r_P = 2$, then $E.\imageofblup{C} = 2$;
but $Q \in E \cap \imageofblup{C}$ and $r_Q \ge 2$, hence the point
$Q$ is the unique singular point of $\imageofblup{C}$ lying over $P$
and $r_{Q} = 2$. Since $\imageofblup{F}.\imageofblup{C} = 1$, it
follows that $Q \not \in \imageofblup{F}$. Therefore $P'=\psi(Q)$ is
a singular point of $C'$ with $r_{P'}=2$ and
\begin{equation*}
\delta_{P'} = \delta_{Q} = \delta_{P} - 1.
\end{equation*}
Note that $E'.C' = 3$ and $r_{P'} = 2$. Therefore there is no
singular points of $C'$ other than $P'$ on $\pi^{-1}(\pi(P)) \cap
C'$.

Suppose that $r_P = 3$, then $E.\imageofblup{C} = 3$; but $Q \in E
\cap \imageofblup{C}$ and $r_Q \ge 2$, hence the point $Q$ is the
unique singular point of $\imageofblup{C}$ lying over $P$. Since
$\imageofblup{F}.\imageofblup{C} = 0$, we have $Q \not \in
\imageofblup{F}$. Therefore $P'$ is a singular point of $C'$ and
\begin{equation*}
\delta_{P'} = \delta_{Q} = \delta_{P} - 3.
\end{equation*}
Furthermore, since $E'.C' = 3$ and $r_{P'} \ge 2$, it follows that
there is no singular points of $C'$ other than $P'$ on
$\pi^{-1}(\pi(P)) \cap C'$.
\end{proof}

The following Lemma is not difficult; but we can't find references.

\begin{lemma}\label{lemma:resolution-of-singularities}
The singularities of a trisection $C$ can be resolved by a finite
sequence of elementary transformations consisting of at most
$\sum_{P \in C} \delta_P$-elementary transformations.
\end{lemma}

\begin{proof}
Let $C$ be a trisection on a ruled surface. Suppose that $P \in X$
is a singular point of $X$. By
Lemma~\ref{lemma:elm_at_singular_point}, we have $\delta_{P'} <
\delta_{P}$. Therefore, after applying finitely many elementary
transformations, we get $\delta_{P'} = 0$, which means $P'$ is a
nonsingular point of $X'$. According to
Lemma~\ref{lemma:elm_at_singular_point}, no singular points other
than infinitely near singular points arise during applying
elementary transformations. Therefore, applying this process to all
singular points of $X$, one can resolve the singularities of $C$.
\end{proof}

\begin{proposition}\label{proposition:general-form-of-bpf-pencils}
Suppose that there exists a nontrivial morphism $h : X \to \PP^1$.
Then
\begin{equation*}
h = p_2 \circ \elm|_X,
\end{equation*}
where $p_2 : Y \times \PP^1 \to \PP^1$ is the second projection and
$\elm$ is a finite sequence of elementary transformations which
transforms $\ruled{E}$ to $Y \times \PP^1$.
\end{proposition}

\begin{proof}
Set $X' = (f \times h)(X) \subset Y \times \PP^1$. The image $X'$ is
a trisection of $Y \times \PP^1$. According to
Lemma~\ref{lemma:resolution-of-singularities}, the singularities of
$X'$ can be resolved by a finite sequence of elementary
transformations. Let $u$ denote a sequence of elementary
transformations which resolves the singularities of $X'$. Set $W =
u(Y \times \PP^1)$ and let $\pi : W \to Y$ be the projection. Note
that $\pi \circ u = p_1$, where $p_1 : Y \times \PP^1 \to Y$ be the
first projection; hence it follows that
\begin{equation*}
\pi \circ u \circ (f \times h) = f.
\end{equation*}
Since the morphism $h$ does not factor through the triple covering
$f$, the image $X'$ is birational to $X$. Therefore $u(X')$ is
isomorphic to $X$. Therefore there is an embedding
\begin{equation*}
j = u \circ (f \times h) : X \hookrightarrow W
\end{equation*}
such that $f = \pi \circ j$. By
Lemma~\ref{preliminaries-lemma:uniqueness-of-embedding}, such an
embedding is unique. Therefore it follows that $W \cong \ruled{E}$
and hence $h = p_2 \circ u^{-1}|_{X}$. Note that $u^{-1}$ is also a
finite sequence of elementary transformations which transforms
$\ruled{E}$ to $Y \times \PP^1$.
\end{proof}

We need the following lemma several times.

\begin{lemma}[{\cite[Lemma~7]{Wolfgang-1992}}]\label{lemma:e-invariant-under-elm}
Let $W$ be a ruled surface with $e$-invariant $e$ and let $P \in W$.
Let $e'$ be the $e$-invariant of $W' = \elm_{P}(W)$. If $P$ is
contained in a minimal degree section of $W$ then $e' = e+1$, and,
if $P$ is not contained any minimal degree sections of $W$, then $e'
= e-1$.
\end{lemma}

We now improve Castelnuovo-Severi inequality for triple coverings.

\begin{theorem}\label{theorem:CS-inequality-for-triple-coverings}
Let $X$ be a smooth irreducible curve of genus $g_x$. Suppose that
there exists a triple covering $f : X \to Y$ where $Y$ is a smooth
irreducible curve of genus $g_y$. Let $\sheaf{E}^{\spcheck}$ be the
Tschirnhausen module for the triple covering $f$ and $e$ the
$e$-invariant of the ruled surface $\ruled{E}$. If there exists a
nontrivial morphism $h : X \to \PP^1$ of degree $d$, then
\begin{equation*}
d \ge \frac{g_x - 3g_y + 2}{2} + \frac{\abs{e}}{2}.
\end{equation*}
\end{theorem}

\begin{proof}
Set $X' = (f \times h)(X) \subset Y \times \PP^1$. Let $P$ be a
singular point of $X'$. According to
Lemma~\ref{lemma:resolution-of-singularities}, there exists a finite
sequence of elementary transformations which resolves the
singularities of $X'$. Set
\begin{equation*}
\nelmresol = \min \{l : \text{$\elm_{P_1, \dotsc, P_l}$ is a
resolution of singularities of $X'$ }\}.
\end{equation*}
Let $\elm_{P_1, \dotsc, P_{\nelmresol}}$ be a minimal sequence of
elementary transformations that is a resolution of singularities of
$X'$. By Lemma~\ref{lemma:resolution-of-singularities}, we have
\begin{equation}\label{equation:l<=sum-of-delta-P}
\nelmresol \le \sum_{P \in X'}{\delta_{P}}.
\end{equation}

We now prove that $\abs{e} \le \nelmresol$. Since $X$ is a
normalization of $X'$, we have an isomorphism $\phi : X \to
\elm_{P_1, \dotsc, P_{\nelmresol}}(X')$. Therefore there exists an
embedding
\begin{equation*}
j = j' \circ \phi : X \hookrightarrow \elm_{P_1, \dotsc,
P_{\nelmresol}}(Y \times \PP^1)
\end{equation*}
such that $f = \pi \circ j$, where $j' : \elm_{P_1, \dotsc,
P_{\nelmresol}}(X') \hookrightarrow \elm_{P_1, \dotsc,
P_{\nelmresol}}(Y \times \PP^1)$ is the inclusion and $\pi :
\elm_{P_1, \dotsc, P_{\nelmresol}}(Y \times \PP^1) \to Y$ is the
projection. However such an embedding $j$ is unique by
Lemma~\ref{preliminaries-lemma:uniqueness-of-embedding}; hence it
follows that $\elm_{P_1, \dotsc, P_{\nelmresol}}(Y \times \PP^1)
\cong \ruled{E}$. The $e$-invariant of $\ruled{E}$ is equal to $e$,
but the $e$-invariant of $Y \times \PP^1$ is equal to zero.
Therefore, by Lemma~\ref{lemma:e-invariant-under-elm}, it follows
that
\begin{equation}\label{equation:|e|<=l}
\abs{e} \le \nelmresol.
\end{equation}
By \eqref{equation:l<=sum-of-delta-P} and
\eqref{equation:|e|<=l}, we have $\abs{e} \le \sum_{P \in
X'}{\delta_P}$. Hence it follows by
\eqref{equation:p_a(X')=2d+3g_y-2} that
\begin{equation*}
2d = g_x - 3g_y +2 + \sum_{P \in X'}{\delta_P} \ge g_x - 3g_y +2 +
\abs{e}. \qedhere
\end{equation*}
\end{proof}

\begin{remarks}\label{remark:final-remarks}\hfill
\begin{enumerate}[(a)]
\item Theorem~\ref{theorem:d>=b_2+4g_y} and
\ref{theorem:CS-inequality-for-triple-coverings} give a partial
answer to Question~\ref{question:converse-of-CS-inequality} in case
of triple covers. On the other hand, as we already remarked in the
introduction, Question~\ref{question:converse-of-CS-inequality} has
been completely answered in case of trigonal curves by the following
results of Maroni.

\begin{maroni}[{\cite{Maroni-1946}}, cf.~{\cite[Proposition~1, Corollary~1]{Martens-Schreyer-1986}}]\label{Maroni's-Result} %
Let $X$ be a trigonal curve of genus $g_x > 4$ with the triple
covering $f : X \to \PP^1$ and let $m_0$ be the so-called Maroni
invariant of $X$. For every integer $d$ with $d \ge g_x - m_0$,
there exists a nontrivial morphism $h : X \to \PP^1$ of degree $d$.
Conversely, if there exists a nontrivial morphism $h : X \to \PP^1$
of degree $d$, then $d \ge g_x - m_0$.
\end{maroni}

\item If $X$ is a trigonal curve, then it is not difficult to prove that
the M-invariant $m$ of $X$ is equal to the Maroni invariant of $X$.
Therefore Theorem~\ref{theorem:d>=b_2+4g_y} and
\ref{theorem:CS-inequality-for-triple-coverings} may be regarded as
a generalization of the results of Maroni to the case of triple
covers of smooth irrational curves.
\end{enumerate}
\end{remarks}

\section{Nontrivial morphisms of minimal possible degree}\label{section:ruled-surface-under-elm}%

We investigate what conditions may be imposed on
$\sheaf{E}^{\spcheck}$ if $X$ admits a nontrivial morphism $X \to
\PP^1$ of minimal possible degree $\frac{g_x - 3g_y + 2}{2} +
\frac{\abs{e}}{2}$;
Proposition~\ref{proposition:minimal-possible=>decomposible}. As a
corollary, we will show that certain triple covers do not admit
nontrivial morphisms of minimal possible degree;
Corollary~\ref{corollary:indecomposable=>not-sharp}. We need the
following proposition.

\begin{proposition}[{cf.~\cite{Keem-Ohbuchi-2004} or \cite{Fuentes-Pedreira-2005}}]\label{proposition:elmof-decomposable} %
Let $\sheaf{F} = \sheaf{O_Y(B_1)} \oplus \sheaf{O_Y}(B_2)$ be a
decomposable rank two vector bundle on a curve $Y$. Let
$S_{\sheaf{F}}$ be a minimal degree section of $\ruled{F}$ and let
$S$ be a section of $\ruled{F}$ such that $S \cap S_{\sheaf{F}} =
\varnothing$. Let $P \in S \cup S_{\sheaf{F}}$. Set $\ruled{G} =
\elm_{P}(\ruled{F})$.
\begin{enumerate}[(a)]
\item If $P \in S$ and $\deg(-B_2 + B_1 + p) \le 0$, then
$\elm_{P}(S_{\sheaf{F}})$ is a minimal degree section of $\ruled{G}$
and
\begin{equation*}
\ruled{G} \cong \PP(\sheaf{O_Y}(B_1) \oplus \sheaf{O_Y}(B_2 - p))
\cong \PP(\sheaf{O_Y}(B_1 + p) \oplus \sheaf{O_Y}(B_2)).
\end{equation*}

\item If $P \in S_{\sheaf{F}}$, then $\elm_{P}(S_{\sheaf{F}})$ is a minimal
degree section of $\ruled{G}$ and
\begin{equation*}
\ruled{G} \cong \PP(\sheaf{O_Y}(B_1 - p) \oplus \sheaf{O_Y}(B_2))
\cong \PP(\sheaf{O_Y}(B_1) \oplus \sheaf{O_Y}(B_2 + p)).
\end{equation*}
\end{enumerate}
\end{proposition}

\begin{proposition}\label{proposition:minimal-possible=>decomposible}
Assume that $e \ge 0$. Suppose that there exists a nontrivial
morphism $h : X \to \PP^1$ of the minimal possible degree $\frac{g_x
- 3g_y + 2}{2} + \frac{e}{2}$. Then the Tschirnhausen module
$\sheaf{E}^{\spcheck}$ is decomposable.
\end{proposition}

\begin{proof}
Set $X' = (f \times h)(X) \subset Y \times \PP^1$ and set $d_0 =
\frac{g_x - 3g_y + 2}{2} + \frac{e}{2}$. It follows by
Remark~\ref{remark:X'-is-a-trisection} that
\begin{equation}\label{equation:sum-of-delta_p<=e}
\sum_{P \in X} \delta_{P} = e.
\end{equation}

\setcounter{step}{0}
\begin{step}\label{step:nelmresol=e}
There exists a sequence $\elm_{P_1, \dotsc, P_e} : Y \times \PP^1
\dashrightarrow \ruled{E}$ consisting of exactly $e$ elementary
transformations that is a resolution of singularities of $X'$.
\end{step}

\begin{proof}[Proof of Step~\ref{step:nelmresol=e}]
According to Lemma~\ref{lemma:resolution-of-singularities}, the
singularities of $X'$ can be resolved by a finite sequence of
elementary transformations with the singular points as centers. Set
\begin{equation*}
\nelmresol = \min \{l : \text{$\elm_{P_1, \dotsc, P_l}$ resolves the
singularities of $X'$} \}.
\end{equation*}
We will prove that $\nelmresol = e$. By
Lemma~\ref{lemma:resolution-of-singularities}, we have
\begin{equation}\label{equation:nelmresol}
\nelmresol \le \sum_{P \in X'} \delta_{P} = e.
\end{equation}

According to \cite[V,~Ex.5.5]{Hartshorne-1977}, there is a finite
sequence of elementary transformations which transform $Y \times
\PP^1$ into $\ruled{E}$; set
\begin{equation*}
\nelmtrans = \min \{ l : \elm_{P_1, \dotsc, P_l}(Y \times \PP^1) =
\ruled{E} \}.
\end{equation*}
Since the $e$-invariants of $Y \times \PP^1$ and $\ruled{E}$ are $0$
and $e$, respectively, it follows by
Lemma~\ref{lemma:e-invariant-under-elm} that
\begin{equation}\label{equation:nelmtrans}
e \le \nelmtrans.
\end{equation}
Let $\elm$ be a sequence of elementary transformations that resolves
the singularities of the image $X'$. Then $\elm(X') \cong X$; hence
\begin{equation}\label{equation:elm(Y*PP^1)=ruled(E)}
\elm(Y \times \PP^1) \cong \ruled{E}
\end{equation}
by Lemma~\ref{preliminaries-lemma:uniqueness-of-embedding}.
Therefore we have
\begin{equation}\label{equation:nelmresol<=nelmtrans}
\nelmtrans \le \nelmresol .
\end{equation}
By \eqref{equation:nelmresol}, \eqref{equation:nelmtrans},
and \eqref{equation:nelmresol<=nelmtrans}, it follows that
\begin{equation}\label{equation:nelmresol=nelmtrans}
\nelmresol = \nelmtrans = e.
\end{equation}
Therefore there exists a sequence of elementary transformations
consisting of exactly $e$ elementary transformations which resolves
the singularities of $X'$.
\end{proof}

\begin{step}\label{step:P_i-contained-in-a-section}
For some $t \in \PP^1$, we have
\begin{equation*}
P_1 \in Y \times \{t\}, \quad P_{i+1} \in \elm_{P_1, \dotsc, P_i}(Y
\times \{ t \})
\end{equation*}
for all $i=1, \dotsc, e-1$.
\end{step}

\begin{proof}[Proof of Step~\ref{step:P_i-contained-in-a-section}]
Let $\elm_{P_1, \dotsc, P_e}$ be a minimal sequence of elementary
transformations that resolves the singularities of $X'$, which
consists of $e$-elementary transformations. Set
\begin{equation*}
X_i = \elm_{P_1, \dotsc, P_i}(X'), \quad W_i = \elm_{P_1, \dotsc,
P_i}(Y \times \PP^1), \quad X_0 = X', \quad W_0 = Y \times \PP^1.
\end{equation*}
Let $e_i$ be the $e$-invariant of $\elm_{P_1, \dotsc, P_i}(Y \times
\PP^1)$. It follows by Lemma~\ref{lemma:e-invariant-under-elm} that
$e_i \le i$, where the equality holds if and only if $P_j$ is
contained in a minimal degree section of $W_{j-1}$ for all $j$ with
$1 \le j \le i$. Note that $\elm_{P_1, \dotsc, P_e}$ consists of
exactly $e$ elementary transformations and the $e$-invariant of
$\ruled{E}$ is equal to $e$; but we have
\begin{equation*}
\elm_{P_1, \dotsc, P_e}(Y \times \PP^1) \cong \ruled{E}
\end{equation*}
by \eqref{equation:elm(Y*PP^1)=ruled(E)}. Therefore it
follows that $ e_i = i$ for all $i$; hence $P_{i+1}$ is contained in
a minimal degree section of $W_i$ for all $i=0, \dotsc, e-1$.

Suppose that $P_1 \in Y \times \{ t \} \subset Y \times \PP^1$ for
some $t \in \PP^1$. Since $Y \times \{ t \}$ is a minimal degree
section of $Y \times \PP^1$, it follows that the section
$\elm_{P_1}(Y \times \{ t \})$ is also a minimal degree section of
$W_1$ by Proposition~\ref{proposition:elmof-decomposable}. However,
since $P_1$ is contained in the minimal degree section $Y \times
\{t\}$, the ruled surface $W_1$ is decomposable by
Proposition~\ref{proposition:elmof-decomposable}. Therefore there
exists a unique minimal degree section of $W_1$ by
\cite[Corollary~1.17]{Maruyama-1970}, which is equal to
$\elm_{P_1}(Y \times \{ t \})$; hence we have $P_2 \in \elm_{P_1}(Y
\times \{ t \})$. Repeating this process, we have $P_{i+1} \in
\elm_{P_1, \dotsc, P_i}(Y \times \{ t \})$ for all $i = 1, \dotsc,
e-1$.
\end{proof}

\noindent \textit{Continue the proof of
Proposition~\ref{proposition:minimal-possible=>decomposible}.} Note
that $\elm_{P_1, \dotsc, P_e}(Y \times \PP^1) \cong \ruled{E}$ by
\eqref{equation:elm(Y*PP^1)=ruled(E)}. By
Proposition~\ref{proposition:elmof-decomposable}, the sections
$\elm_{P_1, \dotsc, P_i}(Y \times \{t \})$ of $\elm_{P_1, \dotsc,
P_i}(Y \times \PP^1)$ are minimal degree sections for all $i$; hence
the ruled surface $\ruled{E}$ is decomposable by
Proposition~\ref{proposition:elmof-decomposable}. Therefore
$\sheaf{E}$ is decomposable.
\end{proof}

\begin{corollary}\label{corollary:indecomposable=>not-sharp}
Assume that $e \ge 0$ and $\sheaf{E}$ is indecomposable. If there
exists a nontrivial morphism $h : X \to \PP^1$ of degree $d$, then
we have
\begin{equation*}
d \gneqq \frac{g_x - 3g_y + 2}{2} + \frac{e}{2}.
\end{equation*}
\end{corollary}

We characterize nontrivial morphisms $X \to \PP^1$ of minimal
possible degree.

\begin{proposition}\label{proposition:triple-cover-admitting-minimal-possible-degree}
Assume that the Tschirnhausen module $\sheaf{E}^{\spcheck}$ is
decomposable; set $\sheaf{E} = \sheaf{O_Y}(B_1) \oplus
\sheaf{O_Y}(B_2)$ with $\deg{B_1} \le \deg{B_2}$. Suppose that there
exists a nontrivial morphism $h : X \to \PP^1$ of degree $\frac{g_x
- 3g_y + 2}{2} + \frac{\abs{e}}{2}$. Then there exist a section $S$
of $\ruled{E}$ which is not a minimal degree section and $e$ points
$Q_1, \dotsc, Q_e \in X \cap S \subset \ruled{E}$ with $\pi_{\ast}
\left( \sum{Q_i} \right) \in \linsys{B_2 - B_1}$ as a divisor on $S$
such that
\begin{equation*}
h = p_2 \circ \elm_{Q_1, \dotsc, Q_e}|_{X},
\end{equation*}
where $p_2 : Y \times \PP^1 \to \PP^1$ is the projection.
\end{proposition}

\begin{proof}
Set $d_0 = \frac{g_x - 3g_y + 2}{2} + \frac{\abs{e}}{2}$. Let $X' =
(f \times h)(X) \subset Y \times \PP^1$. According to
Step~\ref{step:nelmresol=e} of the proof of
Proposition~\ref{proposition:minimal-possible=>decomposible}, there
exists a sequence of elementary transformations
\begin{equation*}
\elm_{P_1, \dotsc, P_e} : Y \times \PP^1 \dashrightarrow \ruled{E}
\end{equation*}
that is a resolution of singularities of $X'$. Set
\begin{equation*}
X_i = \elm_{P_1, \dotsc, P_i}(X'), \quad W_i = \elm_{P_1, \dotsc,
P_i}(Y \times \PP^1), \quad X_0 = X', \quad W_0 = Y \times \PP^1.
\end{equation*}
We may assume that $P_{i+1}$ is a singular point of $X_i$ for all
$i=0, \dotsc, e-1$. Let $r_{i+1}$ be the multiplicity of $X_i$ at
$P_{i+1}$.

First, we will prove that $r_{i+1} = 2$ for all $i$. Suppose
$r_{i+1} = 3$ but $r_{j} = 2$ for all $j \le i$. If we apply
$\elm_{P_{i+1}}$ to $X_i$, then, by
Lemma~\ref{lemma:elm_at_singular_point}, we have
\begin{equation*}
\sum_{P \in X_{i+1}} \delta_{P} = \sum_{P \in X_i} \delta_{P} - 3 =
\sum_{P \in X_0} \delta_{P} - i - 3 = e-i-3.
\end{equation*}
Therefore we need at most ($e-i-3$) elementary transformations to
resolve the singularities of $X_{i+1}$ by
Lemma~\ref{lemma:resolution-of-singularities}. However it
contradicts the assumption: $\elm_{P_1, \dotsc, P_e}$ is the minimal
sequence. Therefore we have $r_{i+1} = 2$ for all $i$.

Set $S_0 = \elm_{P_1, \dotsc, P_e}(Y \times \{ t \})$, which is the
minimal degree section of $\ruled{E}$, and set $Q_{i+1} = \elm_{P_1,
\dotsc, P_e}(F_{i+1})$, where $F_{i+1} = p_1^{-1}(p_1(P_{i+1}))$ is
the fiber of the first projection $p_1 : \elm_{P_1, \dotsc, P_{i}}(Y
\times \PP^1) \to Y$ over $p_1(P_{i+1})$. It is clear that
$\elm_{Q_1, \dotsc, Q_e}$ is the inverse of $\elm_{P_1, \dotsc,
P_e}$; cf.~\cite{Maruyama-1970}. Let $S = \elm_{P_1, \dotsc, P_e}(Y
\times \{ t' \})$ for some $t' \neq t$. Then $S$ is a section of
$\ruled{E}$ such that $S \cap S_0 = \varnothing$ by
Proposition~\ref{proposition:elmof-decomposable}.

Let $\elm_{P_{i+1}} = \bldown_{i+1} \circ \blup_{i+1}^{-1}$. Since
$r_{i+1} = 2$ and $F_{i+1}.X_i = 3$ for all $i$, we have
$\imageofblup{F_{i+1}} \cap \imageofblup{X_i} \neq \varnothing$,
where $\imageofblup{F_{i+1}}$ and $\imageofblup{X_i}$ are strict
transforms under the blowing up $\blup_{i+1}$. Therefore it follows
that
\begin{equation*}
Q_{i+1} \in X \cap \elm_{P_1, \dotsc, P_{i+1}}(Y \times \{ t' \}),
\end{equation*}
hence $Q_i \in X \cap S$ for all $i$. Since $S \cap S_0 =
\varnothing$, it follows by
Proposition~\ref{proposition:elmof-decomposable} that
\begin{equation*}
\elm_{Q_1, \dotsc, Q_e}(\ruled{E}) \cong \PP(\sheaf{O_Y}(B_1 + q_1 +
\dotsb + q_e) \oplus \sheaf{O_Y}(B_2)),
\end{equation*}
where $q_i = \pi(Q_i)$. Since $\elm_{Q_1, \dotsc, Q_e}(\ruled{E})
\cong Y \times \PP^1$, it follows that
\begin{equation*}
\PP(\sheaf{O_Y}(B_1 + q_1 + \dotsb + q_e) \oplus \sheaf{O_Y}(B_2))
\cong \PP(\sheaf{O_Y} \oplus \sheaf{O_Y});
\end{equation*}
hence $\sum q_i \lineqv B_2 - B_1$, that is, $\pi_{\ast} \left(
\sum{Q_i} \right) \lineqv B_2 - B_1$. It is clear that $h = p_2
\circ \elm_{Q_1, \dotsc, Q_e}|_X$.
\end{proof}

\begin{corollary}
Let $Y$ be an irreducible smooth curve of genus $g_y > 1$. For every integer $g_x$ with $g_x \ge 13g_y$ and $e$ with $2g_y - 3 \ge e \ge 0$, there exist an irreducible smooth curve $X$ of genus $g_x$ and a triple covering $f : X \to Y$ such that $X$ does not admit any nontrivial morphisms $X \to \PP^1$ of degree $\frac{g_x - 3g_y + 2}{2} + \frac{\abs{e}}{2}$ and $e$ is the $e$-invariant of the ruled surface $\ruled{E}$ associated to the Tschirnhausen module $\sheaf{E}^{\spcheck}$ for the triple covering $f$.
\end{corollary}

\begin{proof}
Choose an effective divisor $D = Q_1 + \dotsb + Q_e \in
\divisorgroup(Y)$ so that $Q_i \neq Q_j$ for $i \neq j$, $\deg{D} =
e \le 2g_y - 3$, and $H^0(Y, D) = 1$. Choose an effective divisor
$B_1 \in \divisorgroup(Y)$ such that
\begin{equation*}
2\deg{B_1} + \deg{D} = g_x - 3g_y + 2.
\end{equation*}
Set $B_2 = B_1 + D$ and $\sheaf{E} = \sheaf{O_Y}(B_1) \oplus
\sheaf{O_Y}(B_2)$. Let $\pi : \ruled{E} \to Y$ be the projection. We
have
\begin{equation*}
S^3{\sheaf{E}} \otimes \det{\sheaf{E}}^{-1} = \sheaf{O_Y}(2B_1 -
B_2) \oplus \sheaf{O_Y}(B_1) \oplus \sheaf{O_Y}(B_2) \oplus
\sheaf{O_Y}(2B_2 - B_1).
\end{equation*}
By the assumptions $g_x \ge 13g_y$ and $\deg{D} = e \le 2g_y - 3$,
it follows that
\begin{equation*}
\deg(2B_1 - B_2) \ge 2g_y.
\end{equation*}
Hence $S^3{\sheaf{E}} \otimes \det{\sheaf{E}}^{-1}$ is generated by
global sections. Furthermore $H^0(Y, \sheaf{E}^{\spcheck}) = 0$.
According to \cite[Theorem~3.6]{Casnati-Ekedahl-1996}, the zero
locus of a general section contained in $H^0(\ruled{E},
\pi^{\ast}{\det{\sheaf{E}}^{-1}}(3))$ is an irreducible smooth
triple cover of $Y$ with genus $g_x$. Let $S_0$ be the unique
minimal degree section of $\ruled{E}$ and $S$ a section of
$\ruled{E}$ such that $S_0 \cap S = \varnothing$. We may choose a
section $\delta \in H^0(\ruled{E},
\pi^{\ast}{\det{\sheaf{E}}^{-1}}(3))$ so that the zero locus $X$ of
$\delta$ is an irreducible smooth triple cover of $Y$ and $X$ does
not pass through the points in $S \cap \pi^{-1}(\{Q_1, \dotsc,
Q_e\})$.

Let $T$ be a section of $\ruled{E}$ such that $T \cap S_0 =
\varnothing$. Since $S, T \in \linsys{S_0 + \pi^{\ast}{D}}$ by
Lemma~\ref{lemma:two-disjoint-sections}, it follows by
\eqref{preliminaries-equation:pi(S.S)} that
\begin{equation*}
\pi_{\ast}(S \cdot T) \in \linsys{\det{\sheaf{E}} \otimes
\sheaf{L_S}^{-1} \otimes \sheaf{L_{T}}^{-1}} = \linsys{D}.
\end{equation*}
Therefore $S \cap T = \{Q_1, \dotsc, Q_e\}$ because $H^0(Y, D) = 1$.

Suppose that there exists a nontrivial  morphism $h : X \to \PP^1$
of degree $\frac{g_x - 3g_y + 2}{2} + \frac{\abs{e}}{2}$. By
Proposition~\ref{proposition:triple-cover-admitting-minimal-possible-degree},
it follows that $Q_1, \dotsc, Q_e \in X \cap S$, but which
contradicts the choice of $X$. Therefore $X$ does not admit a
nontrivial morphism $h : X \to \PP^1$ of degree $\frac{g_x - 3g_y +
2}{2} + \frac{\abs{e}}{2}$.
\end{proof}

\section{Existence of triple covers}\label{section:existence-of-triple-covers}

We give a partial answer to Question~\ref{question:weakened-version-of-converse-of-CS-inequality}; the existence of triple covers, which admit nontrivial morphisms of certain degrees, of a given base curve. We first investigate behaviors of trisections under elementary transformations. Let $C$ be a trisection on a ruled surface $\ruled{F}$. Assume that $C \lineqv 3S + \pi^{\ast}{Z}$ as a divisor in $\ruled{F}$, where $S$ is a (not necessarily minimal degree) section of $\ruled{F}$, $\pi : \ruled{F} \to Y$ is the projection, and $Z \in \divisorgroup{Y}$. Let $P \in \ruled{F}$. Set
\begin{equation*}
p = \pi(P), \quad S' = \elm_{P}(S), \quad C' = \elm_{P}(C).
\end{equation*}

\begin{proposition}\label{proposition:strict-transformation-of-X}
Assume that $P \in S$. If $P \not \in C$ then $C' \lineqv 3S' +
\pi^{\ast}(Z + 3p)$, and if $P \in C$ then $C' \lineqv 3S' +
\pi^{\ast}(Z + 2p)$. Assume that $P \not \in S$. If $P \not \in C$
then $C' \lineqv 3S' + \pi^{\ast}Z$, and if $P \in C$ then $C'
\lineqv 3S' + \pi^{\ast}(Z - p)$.
\end{proposition}

\begin{proof}
We use similar techniques in \cite{Keem-Ohbuchi-2004} which deals
the behavior of double covers under elementary transformations. Set
$\elm_{P} = \bldown \circ \blup^{-1}$ and
$P'=\bldown(\imageofblup{F})$, where $F = \pi^{-1}(\pi(P))$. Suppose
$P \in S$ and $P \not \in C$. Then
\begin{equation*}
\begin{split}
\bldown^{\ast}C' &= \imageofblup{C} + 3 \imageofblup{F} =
\blup^{\ast}C + 3\imageofblup{F} = 3 \imageofblup{S} + 3E +
\blup^{\ast} \pi^{\ast}Z
+ 3\imageofblup{F} \\ %
&= 3\bldown^{\ast}S' + \bldown^{\ast}\pi^{\ast}(Z + 3p) =
\bldown^{\ast}(3S' + \pi^{\ast}(Z + 3p)).
\end{split}
\end{equation*}

Suppose $P \in S$ and $P$ is a nonsingular point of $C$. Then
\begin{equation*}
\begin{split}
\bldown^{\ast}C' &= \imageofblup{C} + 2 \imageofblup{F} =
\blup^{\ast}C - E + 2\imageofblup{F} = 3 \imageofblup{S} + 3E +
\blup^{\ast} \pi^{\ast}Z
- E + 2\imageofblup{F} \\ %
&= 3\bldown^{\ast}S' + \bldown^{\ast}\pi^{\ast}(Z + 2p) =
\bldown^{\ast}(3S' + \pi^{\ast}(Z + 2p)).
\end{split}
\end{equation*}

Suppose $P \not \in S$ and $P \not \in C$. Then
\begin{equation*}
\begin{split}
\bldown^{\ast}C' &= \imageofblup{C} + 3 \imageofblup{F} =
\blup^{\ast}C + 3\imageofblup{F} = 3 \imageofblup{S} + \blup^{\ast}
\pi^{\ast}Z +
3\imageofblup{F} \\ %
&= 3(\bldown^{\ast}S' - \imageofblup{F}) + \bldown^{\ast}\pi^{\ast}Z
+ 3 \imageofblup{F} = \bldown^{\ast}(3S' + \pi^{\ast}Z).
\end{split}
\end{equation*}

Suppose $P \not \in S$ and $P \in C$. Then
\begin{equation*}
\begin{split}
\bldown^{\ast}C' &= \imageofblup{C} + 2 \imageofblup{F} =
\blup^{\ast}C - E + 2\imageofblup{F} = 3 \imageofblup{S} +
\blup^{\ast} \pi^{\ast}Z - E
+ 2\imageofblup{F} \\ %
&= 3(\bldown^{\ast}S' - \imageofblup{F}) +
\bldown^{\ast}\pi^{\ast}(Z - p + p) - E + 2 \imageofblup{F} =
\bldown^{\ast}(3S' + \pi^{\ast}(Z-p)).
\end{split}
\end{equation*}
\end{proof}

We need the following Lemma.

\begin{lemma}\label{lemma:elementary-transformation-of-sections}
Let $S$ and $T$ be sections of a ruled surface $\ruled{F}$ over $Y$
with $T \lineqv S + \pi^{\ast}{Z}$ for some $Z \in
\divisorgroup(Y)$, where $\pi : \ruled{F} \to Y$ is the projection.
Let $\elm_{P}$ be the elementary transformation with center $P \in
S$ and let $S'$, $T'$ be the image $\elm_{P}(S)$, $\elm_{P}(T)$,
respectively. Then we have $T' \lineqv S' + \pi^{\ast}(Z + \pi(P))$.
\end{lemma}

\begin{proof}
Set $\elm_{P} = \bldown \circ \blup^{-1}$ and
$P'=\bldown(\imageofblup{F})$, where $F = \pi^{-1}(\pi(P))$. Then we
have
\begin{equation*}
\begin{split}
\psi^{\ast}{T'} &= \imageofblup{T} + \imageofblup{F} =
\phi^{\ast}{T} + \imageofblup{F} \lineqv \phi^{\ast}{S} +
\phi^{\ast}(\pi^{\ast}{Z}) + \imageofblup{F} \\ %
&= \imageofblup{S} + E + \phi^{\ast}(\pi^{\ast}{Z}) +
\imageofblup{F} = \psi^{\ast}{S'} + \psi^{\ast}(\pi^{\ast}(Z + p)).
\end{split}
\end{equation*}
\end{proof}

The following theorem provides a partial answer to
Question~\ref{question:weakened-version-of-converse-of-CS-inequality}.

\begin{theorem}\label{theorem:existence-irrational-curve}
Let $Y$ be an irreducible smooth curve of genus $g_y \ge 1$, and let
$g_x$ be an integer with $g_x \ge 37g_y - 2$. For every integer $d$
with
\begin{equation*}
d \ge \frac{g_x - 3g_y + 2}{2} + g_y,
\end{equation*}
there exists an irreducible smooth curve $X$ of genus $g_x$ and a
triple covering $f : X \to Y$ such that $X$ admits a nontrivial
morphism $h : X \to \PP^1$ of degree $d$. If $g_y \ge 5$, we get
sharper lower bound for $d$:
\begin{equation*}
d \ge \frac{g_x - 3g_y + 2}{2} + \frac{g_y + 3}{2}.
\end{equation*}
\end{theorem}

\begin{lemma}\label{lemma:one-dimensional-family-of-sections}
Let $\sheaf{E} = \sheaf{O_Y} \oplus \sheaf{O_Y}(D)$ be a rank $2$
vector bundle on an irreducible smooth curve $Y$ and let $E = y_1 +
\dotsb + y_n \in \linsys{D}$ be an effective divisor consisting of
distinct points. Let $S$ be a section of $\ruled{E}$ with $S \cap
S_0 = \varnothing$, where $S_0$ is a minimal degree section. Set
$P_i = S \cap \pi^{-1}(y_i)$. Then, for the pair $(S, \{ P_1,
\dotsc, P_n\})$, there exists an one-dimensional family $V$ of
sections of $\ruled{E}$ such that, for all $T \in V$,
\begin{equation*}
S \in V, \quad \text{$P_i \in T$ for all $i$}, \quad T \in
\linsys{S_0 + \pi^{\ast}{D}}.
\end{equation*}
\end{lemma}

\begin{proof}
Let $\elm_{P_i} = \bldown_i \circ \phi_i^{-1}$ be the elementary
transformation of $\ruled{E}$ with center $P_i$. Set $Z_i =
\bldown(\imageofblup{F_i})$ where $F_i = \pi^{-1}(P_i)$ and
$\imageofblup{F_i}$ is the strict transform of $F_i$ under the
blowing-up $\blup$. By
Proposition~\ref{proposition:elmof-decomposable}, we have
\begin{equation*}
\elm_{P_1, \dotsc, P_n}(\ruled{E}) = Y \times \PP^1, \quad
\elm_{Z_1, \dotsc, Z_n}(Y \times \PP^1) = \ruled{E}.
\end{equation*}
Note that
\begin{equation*}
\elm_{P_1, \dotsc, P_n}{S}= Y \times \{ a \}, \quad \elm_{P_1,
\dotsc, P_n}{S_0}= Y \times \{ b \}
\end{equation*}
for some $a, b \in \PP^1$ because they are sections of $Y \times
\PP^1$. Set
\begin{equation*}
V = \{ \elm_{Z_1, \dotsc, Z_n}(Y \times \{ c \}) \ | \ c \in \PP^1,
c \neq b \}.
\end{equation*}
It is clear that $S \in V$. Let $T = \elm_{Z_1, \dotsc, Z_n}(Y
\times \{ c \}) \in V$. Since $Z_i \notin Y \times \{ c \}$ for all
$i$, we have $P_i \in T$ for all $i$. By
Lemma~\ref{lemma:elementary-transformation-of-sections}, we have $T
\lineqv S_0 + \pi^{\ast}{E}$. Therefore $V$ is the desired family.
\end{proof}

\begin{remarkof}\label{remark:no-other-intersection}
For any $S_1, S_2 \in V$ we have $S_1 \cap S_2 = \{ P_1, \dotsc, P_n
\}$ since $S_1.S_2 = n$.
\end{remarkof}

\begin{proof}[{Proof of Theorem~\ref{theorem:existence-irrational-curve}}]
Set
\begin{equation*}
d_0 =
\begin{cases}
(g_x - 3g_y + 2)/{2} & \text{if $g_x - 3g_y + 2 \equiv 0 \pmod{2}$,} \\ %
(g_x - 3g_y + 2)/{2} + {1}/{2} & \text{if $g_x - 3g_y + 2 \equiv 1
\pmod{2}$.}
\end{cases}
\end{equation*}
Suppose that $g_x - 3g_y + 2 \equiv 0 \pmod{2}$. Let $t$ be an
integer with
\begin{equation*}
g_y \le t \le \frac{d_0 - 2g_y}{3}.
\end{equation*}
Choose a divisor $D = y_1 + \dotsb + y_{2t}$, $\deg{D} = 2t$ of $Y$
consisting of distinct points. Choose an effective divisor $D'$ of
$Y$ with $\deg{D'} = d_0 - 3t$. Set
\begin{equation*}
B_1 = D + D', \quad B_2 = 2D + D', \quad \sheaf{E} =
\sheaf{O_Y}(B_1) \oplus \sheaf{O_Y}(B_2).
\end{equation*}
Let $\pi : \ruled{E} \to Y$ be the projection. Since $\sheaf{E}$ is
decomposable, there exists a \emph{unique} minimal degree section
$S_0$ of $\ruled{E}$ by \cite[Corollary~1.17]{Maruyama-1970} and
there exists a section $S \in \linsys{S_0 + \pi^{\ast}D}$ with $S
\cap S_0 = \varnothing$ by Lemma~\ref{lemma:two-disjoint-sections}.
Set
\begin{equation*}
S \cap \pi^{-1}(y_i) = \{P_i\}.
\end{equation*}

Since $\deg{D} \ge 2g_y$, there exists a divisor $E = y_1' + \dotsb
+ y_{2t}' \in \linsys{D}$ such that $y_i' \neq y_j'$ for $i \neq j$
and $\supp{E} \cap \supp{D} = \varnothing$. Set
\begin{equation*}
S \cap \pi^{-1}(y_i') = \{ P_i' \}.
\end{equation*}
Let $V$ be the one dimensional family of sections given by
Lemma~\ref{lemma:one-dimensional-family-of-sections} corresponding
to $(E, \{P_1', \dotsc, P_{2t}'\})$. Let $S_1, S_2 \in V$ be
sections with $S_1 \neq S_2$ and $S_i \neq S$ for $i=1,2$. Set
\begin{equation*}
S_1 \cap \pi^{-1}(y_i) = \{ Q_i \}, \quad S_2 \cap \pi^{-1}(y_i) =
\{ R_i \}.
\end{equation*}
By Remark~\ref{remark:no-other-intersection}, we have $Q_i \neq
R_i$; cf.~Figure~\ref{figure:sections}.
\begin{figure}
\centering{\includegraphics[scale=1]{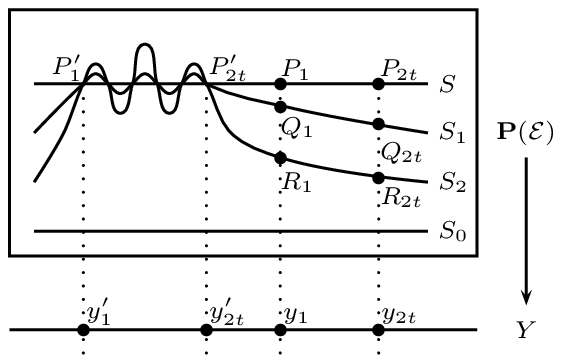}}%
\caption{Sections on $\ruled{E}$}\label{figure:sections}
\end{figure}

Define a linear system $H$ on $\ruled{E}$ by
\begin{equation*}
H = \left\{ X \in \linsys{3S_0 + \pi^{\ast}(2B_2 - B_1)} \ : \
\text{$P_i, Q_i, R_i \in X $ for all $i$} \right\}.
\end{equation*}

\setcounter{step}{0}
\begin{step}\label{step:S_QandS_R}
There exist irreducible smooth curves $S_Q$, $S_R$ such that $S_Q,
S_R \in \linsys{S_0 + \pi^{\ast}B_1}$ and $Q_i \in S_Q$, $R_i \in
S_R$ for all $i=1, \dotsc, 2t$.
\end{step}

\begin{proof}[Proof of Step~\ref{step:S_QandS_R}]
Let $\elm_{Q_i} = \bldown_{i} \circ \blup{_i}^{-1}$ be the
elementary transformation with center $Q_i$. Set $Q_i' =
\bldown_{i}(\imageofblup{F})$, where $\imageofblup{F}$ is the strict
transform by $\blup_{i}$ of the fiber $F = \pi^{-1}(Q_i)$. Set
\begin{equation*}
T_0 = \elm_{Q_1, \dotsc, Q_{2t}}(S_0).
\end{equation*}
Since $S_1$ is a section of $\ruled{E}$ such that $S_1 \cap S_0 =
\varnothing$, it follows by
Proposition~\ref{proposition:strict-transformation-of-X} that
\begin{equation*}
\elm_{Q_1, \dotsc, Q_{2t}}(\ruled{E}) \cong Y \times \PP^1;
\end{equation*}
hence we have $T_0$ is a minimal degree section of $Y \times \PP^1$.
Define a subset $H_Q$ of the linear series $\linsys{T_0 +
\pi^{\ast}(D')}$ on $Y \times \PP^1$ by
\begin{equation*}
H_Q = \{ S_Q' \in \linsys{T_0 + \pi^{\ast}(D')} \ : \ \text{$Q_i'
\notin S_Q'$ for all $i$}\}.
\end{equation*}

Since $\deg{D'} \ge 2g_y + 1$, the linear series $\linsys{S_0 +
\pi^{\ast}(D')}$ is very ample by \cite[V. Ex.
2.11]{Hartshorne-1977}. Therefore there exists an irreducible smooth
curve $S_Q' \in H_Q$. Set
\begin{equation*}
S_Q = \elm_{Q_1', \dotsc, Q_{2t}'}(S_Q')
\end{equation*}
By Lemma~\ref{lemma:elementary-transformation-of-sections}, we have
$S_Q \in \linsys{S_0 + \pi^{\ast}{B_1}}$. Since $Q_i' \notin S_Q'$,
it is clear that $Q_i \in S_Q$ for all $i=1, \dotsc, 2t$. Therefore
we get the desired $S_Q$. Applying the same method to $R_i$, we get
the desired $S_R$.
\end{proof}

\begin{step}\label{step:general-members-are-smooth-triple-coverings}
A general member $X$ of $H$ is an irreducible smooth triple cover of
$Y$ of genus $g_x$.
\end{step}

\begin{proof}[Proof of Step~\ref{step:general-members-are-smooth-triple-coverings}]
Let $V_P$, $V_Q$, $V_R$ be the one dimensional families of sections
corresponding to $(D, \{P_1, \dotsc, P_{2t}\})$, $(D, \{Q_1, \dotsc,
Q_{2t}\})$, $(D, \{R_1, \dotsc, R_{2t}\})$ given by
Lemma~\ref{lemma:one-dimensional-family-of-sections}, respectively.
Then, for any $S' \in V$, $S_1 \in V_Q$, $S_2 \in V_R$, we have
\begin{equation*}
S' + S_1 + S_R \in H \text{ and } S' + S_Q + S_2 \in H.
\end{equation*}
By Remark~\ref{remark:no-other-intersection}, it follows that
\begin{align*}
\bp(V) &= \{P_1, \dotsc, P_{2t}\}, \\ %
\bp(V_Q) &= \{Q_1, \dotsc, Q_{2t}\}, \\ %
\bp(V_R) &= \{R_1, \dotsc, R_{2t}\}.
\end{align*}
Therefore we have
\begin{equation*}
\bp(H) = \{P_1, \dotsc, P_{2t}, Q_1, \dotsc, Q_{2t}, R_1, \dotsc,
R_{2t}\}.
\end{equation*}
By Bertini theorem of characteristic zero (cf.~\cite[III,
10.9]{Hartshorne-1977}), general members of $H$ is smooth outside
$\bp(H)$.

Since $S' + S_1 + S_R$ and $S' + S_Q + S_2$ do not contain any fiber
of $\pi$, they are singular triple coverings of $Y$. Note that
$$S^3\sheaf{E} \otimes \det{\sheaf{E}}^{-1} = \sheaf{O_Y}(2B_2-B_1) \oplus \sheaf{O_Y}(B_2) \oplus \sheaf{O_Y}(B_1) \oplus \sheaf{O_Y}(2B_1-B_2),$$
where
\begin{equation*}
\deg(2B_2 - B_1), \deg(B_2), \deg(B_1), \deg(2B_1-B_2) > 2g_y - 1
\end{equation*}
by the assumptions $g_x \ge 37g_y - 2$ and the choice of $d'$: $d'
\ge 2g_y$. Therefore $S^3\sheaf{E} \otimes \det{\sheaf{E}}^{-1}$ is
globally generated; hence, by
\cite[Thoerem~3.6]{Casnati-Ekedahl-1996}, general members of $H$ are
(possibly singular) triple coverings of $Y$. Since $h^0(Y,
\sheaf{E}^{\spcheck}) = 0$, it follows by
\cite[Theorem~3.6]{Casnati-Ekedahl-1996} that every triple covers in
$H$ are connected. Therefore we proved that general members of $H$
are connected triple covers of $Y$ which are smooth outside
$\bp(H)$.

Let $X$ be a general member of $H$ which is smooth outside $\bp(H)$.
Since
\begin{equation*}
X \cap \pi^{-1}(y_i) = \{ P_i, Q_i, R_i \},
\end{equation*}
it follows that $X$ cannot have a singular point on $P_i$, $Q_i$,
and $R_i$; hence $X$ is smooth. Thus a general member of $H$ is an
irreducible smooth triple cover of $Y$. Let $X \in H$ be an
irreducible \emph{smooth} triple covering of $Y$ with the triple
covering map $f : X \to Y$. By \cite[Proposition~8.1]{Miranda-1985},
the vector bundle $\sheaf{E}^{\spcheck}$ is the Tschirnhausen module
for $f : X \to Y$. Therefore it follows by
\cite[Proposition~4.7]{Miranda-1985} and Riemann-Hurwitz formula
that the genus of $X$ is equal to $g_x$.
\end{proof}

\begin{step}\label{step:higher-degree}
Let $X \in H$ be a smooth triple cover of $Y$ with the triple
covering $f : X \to Y$. For every integer $d$ with
\begin{equation*}
d = d_0 + t, \quad d \ge d_0 + t + 4g_y,
\end{equation*}
there exists a nontrivial morphism $h : X \to \PP^1$ of degree $d$.
\end{step}

\begin{proof}[Proof of Step~\ref{step:higher-degree}]
Set
\begin{equation*}
T = P_1 + \dotsb + P_{2t} \in \divisorgroup(S).
\end{equation*}
Choose an effective divisor $T_1$ on $S$ satisfying the followings:
$T_1 \le S \cap X$, and if $aP \le T_1$ and $bP \le T$ for some $P
\in S$ and $a, b \ge 0$, then $(a+b)P \le S \cap \pi^{\ast}{E}$.
Choose an effective divisor $T_2$ on $S$ satisfying
$$\supp{T_2} \subset S - S \cap X \text{ and } \deg{T_1} + \deg{T_2} \ge
2g_y.$$
Choose an effective divisor $T_3$ on $S_0$ satisfying
\begin{equation*}
\pi_{\ast}{T_3} \lineqv \pi_{\ast}(T_1 + T_2), \quad \supp{T_3} \cap
\supp(T_1 + T_2) = \varnothing,
\end{equation*}
which is possible because $\deg{T_1} + \deg{T_2} \ge 2g_y$. Set $t_i
= \deg{T_i}$($i = 1, 2$).

First of all, consider $\elm_{T}$. By
Proposition~\ref{proposition:elmof-decomposable},
$\elm_{T}(\ruled{E}) \cong Y \times \PP^1$. Set $X' = \elm_{T}{X}$.
By Proposition~\ref{proposition:strict-transformation-of-X},
\begin{equation*}
X' \lineqv 3Y_0 + \pi^{\ast}(2B_2 - B_1 - D) \numeqv 3Y_0 + (d_0 +
t)F.
\end{equation*}
Set $h = p_2 \circ \elm_{T}|_{X}$. Then $h$ is a nontrivial morphism
of degree $d_0 + t$.

Consider now $\elm_{T_3, T_1 + T_2, T}$. By
Proposition~\ref{proposition:elmof-decomposable}, it follows that
\begin{equation*}
\elm_{T_3, T_1 + T_2, T}{\ruled{E}} \cong Y \times \PP^1.
\end{equation*}
Set $X' = \elm_{T_3, T_1 + T_2, T}(X)$. By
Proposition~\ref{proposition:strict-transformation-of-X},
\begin{equation*}
X' \lineqv 3Y_0 + \pi^{\ast}(2B_2 - B_1 - T - T_1 + 3T_3) \numeqv
3Y_0 + (d_0 + t + 2t_1 + 3t_2)F.
\end{equation*}
Set $h = p_2 \circ \elm_{T}|_{X}$. Then $h$ does not factor through
the triple covering $f$ and we have
\begin{equation*}
\deg{h} = d_0 + t + 2t_1 + 3t_2.
\end{equation*}
Note that, by the assumption $g_x \ge 37g_y - 2$, we have
\begin{equation*}
2g_y \le t_1 \le 2b_2 - b_1.
\end{equation*}
We may choose $T_2$ such that $\deg{T_2}$ is arbitrary large.
Therefore any integer greater than or equal to $4g_y$ can be
represented by $2t_1 + 3t_2$ for some $t_1$ and $t_2$.

Note that we have $\frac{d_0 - 2g_y}{3} \ge 5g_y$ by the assumption
$g_x \ge 37g_y - 2$. Therefore it follows that
\begin{equation*}
\{ d : d \ge d_0 + g_y \} = \left\{ d_0 + t : g_y \le t \le
\frac{d_0 - 2g_y}{3} \right\} \cup \{d : d \ge d_0 + 5g_y \}.
\end{equation*}

Suppose that $g_x - 3g_y + 2 \equiv 1 \pmod{2}$. Let $t$ be an
integer with $g_y + 1 \le t \le \frac{d_0 - 2g_y}{3}$ and fix a
divisor $D = y_1 + \dotsb + y_{2t-1}$ of $Y$ with $\deg{D} = 2t-1$
which consists of distinct points and then repeat the above proof.
\end{proof}

\noindent \textit{Continue the proof of
Theorem~\ref{theorem:existence-irrational-curve}.} Hence far, we
proved the existence of triple coverings and the existence of
nontrivial morphisms for $d \ge \frac{g_x - 3g_y + 2}{2} + g_y$. It
remains to prove the existence for $d \ge \frac{g_x - 3g_y + 2}{2} +
\frac{g_y+3}{2}$ if $g_y \ge 5$. Suppose that $g_y \ge 2$. A
classical theorem of Halphen says that a curve $Y$ of genus $g_y \ge
2$ has a nonspecial very ample divisor $D$ with $\deg{D} \ge g_y +
3$. Therefore, in the above proof, we can take $t$ with $t \ge
\frac{g_y + 3}{2}$ and take this divisor $D$ instead of arbitrary
divisor of degree $2t$. Then we get better lower bound $d \ge d_0 +
\frac{g_y + 3}{2}$.
\end{proof}

\subsection*{Acknowledgements}

This paper covers a part of the content of Ph.D.~thesis of the author submitted to Seoul National University, Korea. The author would like to thanks to Prof.~Changho Keem, his advisor, for valuable comments and discussions. The author is also grateful to the referee for pointing out some inaccuracies in the previous version and the very interesting Remark~\ref{remark:referee}.

\providecommand{\bysame}{\leavevmode\hbox to3em{\hrulefill}\thinspace}
\providecommand{\MR}{\relax\ifhmode\unskip\space\fi MR }
\providecommand{\MRhref}[2]{%
  \href{http://www.ams.org/mathscinet-getitem?mr=#1}{#2}
}
\providecommand{\href}[2]{#2}


\begin{thebibliography}{10}

\bibitem{ACGH-1985}
E.~Arbarello, M.~Cornalba, P.~Griffiths, and J.~Harris, \emph{Geometry of
  algebraic curves. {V}ol. {I}}, Springer-Verlag, 1985. \MR{86h:14019}

\bibitem{Ballico-Keem-Park-2004}
E.~Ballico, C.~Keem, and S.~Park, \emph{Double covering of curves}, Proc. Amer.
  Math. Soc. \textbf{132} (2004), no.~11, 3153--3158 (electronic).
  \MR{2005a:14037}

\bibitem{Casnati-Ekedahl-1996}
G.~Casnati and T.~Ekedahl, \emph{Covers of algebraic varieties. {I}. {A}
  general structure theorem, covers of degree {$3,4$} and {E}nriques surfaces},
  J. Algebraic Geom. \textbf{5} (1996), no.~3, 439--460. \MR{97c:14014}

\bibitem{Coppens-Keem-Martens-1992}
M.~Coppens, C.~Keem, and G.~Martens, \emph{Primitive linear series on curves},
  Manuscripta Math. \textbf{77} (1992), no.~2-3, 237--264. \MR{93j:14028}

\bibitem{Fuentes-Pedreira-2005}
L.~Fuentes and M.~Pedreira, \emph{The projective theory of ruled surfaces},
  Note Mat. \textbf{24} (2005), no.~1, 25--63. \MR{2199622}

\bibitem{Hartshorne-1977}
R.~Hartshorne, \emph{Algebraic geometry}, Springer-Verlag, 1977, Graduate Texts
  in Mathematics, No. 52. \MR{57 \#3116}

\bibitem{Kato-Keem-Ohbuchi-1997}
T.~Kato, C.~Keem, and A.~Ohbuchi, \emph{On triple coverings of irrational
  curves}, Tsukuba J. Math. \textbf{21} (1997), no.~2, 421--441. \MR{98m:14031}

\bibitem{Keem-Ohbuchi-2004}
C.~Keem and A.~Ohbuchi, \emph{On the {C}astelnuovo-{S}everi inequality for a
  double covering}, preprint, 2004.

\bibitem{Maroni-1946}
A.~Maroni, \emph{Le serie lineari speciali sulle curve trigonali}, Ann. Mat.
  Pura Appl. (4) \textbf{25} (1946), 343--354. \MR{0024182 (9,463j)}

\bibitem{Martens-Schreyer-1986}
G.~Martens and F.-O. Schreyer, \emph{Line bundles and syzygies of trigonal
  curves}, Abh. Math. Sem. Univ. Hamburg \textbf{56} (1986), 169--189.
  \MR{88d:14019}

\bibitem{Maruyama-1970}
M.~Maruyama, \emph{On classification of ruled surfaces}, Lectures in
  Mathematics, Department of Mathematics, Kyoto University, vol.~3, Kinokuniya
  Book-Store Co. Ltd., Tokyo, 1970. \MR{43 \#1990}

\bibitem{Miranda-1985}
R.~Miranda, \emph{Triple covers in algebraic geometry}, Amer. J. Math.
  \textbf{107} (1985), no.~5, 1123--1158. \MR{86k:14008}

\bibitem{Wolfgang-1992}
W.~Seiler, \emph{Deformations of ruled surfaces}, J. Reine Angew. Math.
  \textbf{426} (1992), 203--219. \MR{93b:14063}
\end{thebibliography}
\end{document}